\documentclass[10pt]{amsart}
\usepackage[cp1251]{inputenc}
\usepackage[english]{babel}
\usepackage{amsmath}
\usepackage{amssymb}
\usepackage{amsfonts}

\begin{document}
\renewcommand{\refname}{References}

\thispagestyle{empty}

\title[Numerical Simulation of 2.5-Set of Iterated Ito Stochastic Integrals]
{Numerical Simulation of 2.5-Set of Iterated Ito Stochastic 
Integrals of Multiplicities 1 to 5 From the Taylor--Ito Expansion}
\author[D.F. Kuznetsov]{Dmitriy F. Kuznetsov}
\address{Dmitriy Feliksovich Kuznetsov
\newline\hphantom{iii} Peter the Great Saint-Petersburg Polytechnic University,
\newline\hphantom{iii} Polytechnicheskaya ul., 29,
\newline\hphantom{iii} 195251, Saint-Petersburg, Russia}%
\email{sde\_kuznetsov@inbox.ru}

\vspace{5mm}

\maketitle {\small
\begin{quote}
\noindent{\sc Abstract.} 
The article is devoted to the construction of effective procedures of
the mean-square approximation of iterated Ito
stochastic integrals of multiplicities 1 to 5 from the
Taylor--Ito expansion
based on multiple Fourier--Legendre series.  
The results of the article can be applied to
the 
implementation of numerical methods with 
orders 1.5,
2.0, and 2.5 of strong convergence for Ito stochastic differential equations  
with multidimensional non-commutative noise.
\medskip
\end{quote}
}

\setlength{\baselineskip}{2.0em}

\tableofcontents

\setlength{\baselineskip}{1.2em}

\section{Introduction}

\vspace{5mm}

This paper is a continuation of the author's research \cite{18a},
\cite{arxiv-14}
on numerical methods with rather high
orders of strong convergence for Ito stochastic differential 
equations (SDEs). The definition of
strong convergence will be given below.
The development of such numerical methods is topical due to 
a wide range of applications for
Ito SDEs \cite{KlPl2}-\cite{KPS}. 
In particular, these equations arise in optimal 
stochastic control, signal
filtering against the background of random noises, 
parameter estimation for stochastic systems as
well as in stochastic stability and bifurcations analysis 
\cite{KlPl2}-\cite{KPS}. 
Also the Ito SDEs represent adequate
mathematical models for dynamic systems of different physical origin 
that are affected by random
perturbations. They are used as mathematical models in stochastic 
mathematical finance,
hydrology and seismology, geophysics, chemical kinetics 
and population dynamics,
electrodynamics, medicine and other fields (see \cite{KlPl2}-\cite{KPS}).
On the other hand, new numerical methods with rather high orders 
of strong convergence are
needed for the Ito SDEs because one of the elementary numerical 
methods --- the Euler scheme has
insufficient accuracy for a series of practical problems under 
standard assumptions; the details can
be found in \cite{KlPl2}.

This paper follows a promising approach to the numerical 
integration of Ito SDEs \cite{KlPl2}, \cite{Mi2}-\cite{KPS}
that is based on the stochastic analogs of the Taylor formula 
(the so-called Taylor--Ito and Taylor--Stratonovich expansions
\cite{KlPl2}, \cite{Mi2}-\cite{2006}) 
for the solutions of Ito SDEs. This approach includes the finite
partitioning of the time variable and also the numerical solution 
of an Ito SDE at discrete moments of time using
the stochastic analogs of the Taylor formula. The numerical methods
with orders 1.5, 2.0, and 2.5 of strong convergence 
will be considered in the article.

This paper employs the so-called unified Taylor--Ito expansion 
\cite{kk5}-\cite{2018aa1} with a minimum set of iterated
Ito stochastic integrals, which is a simplifying factor at 
the implementation stage of numerical
methods. The iterated Ito stochastic integrals figuring in the 
numerical schemes with a strong convergence
of orders 1.5, 2.0, and 2.5 are approximated using the method of 
generalized multiple 
Fourier series, which
was considered in a series of papers of the author 
\cite{2006}-\cite{new-art-1-xxy}. 
As was noted in \cite{18a},
\cite{arxiv-14}, this method
does not require the splitting of an integration interval $[t, T]$ of the 
iterated Ito stochastic integrals;
recall that its length $T-t$ gives the integration step of the 
numerical methods for Ito SDEs and
hence is a sufficiently small value. In accordance with experimental 
results \cite{2006}, the splitting of
the interval $[t, T]$ leads to an inadmissibly high computational cost. 
As a rule, this splitting is used
in the approximation methods of iterated Ito stochastic integrals
based on integral sums \cite{Mi2}, \cite{Mi3}, \cite{allen}.

As mentioned in \cite{18a},
\cite{arxiv-14}, in a number of publications \cite{KlPl2}, \cite{Mi2}, 
\cite{Mi3},
numerical schemes
with strong convergence of high orders (1.5, 2.0, and 2.5) for 
the Ito SDEs have been proposed.
However, these methods do not contain
efficient
mean-square approximation procedures for the iterated Ito stochastic 
integrals in
the case of multidimensional non-commutative noises. Generally, the 
authors \cite{KlPl2}, \cite{Mi2}, \cite{Mi3}
introduced
some simplifying assumptions on the additivity, commutativity or 
smallness of the noises, which
results in a considerable simplification of the numerical 
modeling problem of iterated Ito stochastic
integrals. Like \cite{18a}, \cite{arxiv-14}, 
this paper will partially eliminate this drawback.

Note that the properties of numerical schemes for the Ito SDEs 
(including the ones with a
strong convergence of orders 1.5, 2.0, and 2.5) 
were well studied in \cite{KlPl2}, \cite{Mi2}, \cite{Mi3}.
In particular, their stability was
thoroughly analyzed. The goal of this paper is to develop 
efficient numerical modeling procedures
for the iterated Ito stochastic integrals of multiplicities 
1 to 5 including the exact calculation and
efficient estimation of the mean-square approximation errors of 
these stochastic integrals.

Let $(\Omega,$ ${\rm F},$ ${\sf P})$ be a complete probability space, let 
$\{{\rm F}_t, t\in[0,T]\}$ be a nondecreasing 
right-continous family of $\sigma$-subfields of ${\rm F},$
and let ${\bf f}_t$ be a standard $m$-dimensional Wiener 
stochastic process, which is
${\rm F}_t$-measurable for any $t\in[0, T].$ We assume that the components
${\bf f}_{t}^{(i)}$ $(i=1,\ldots,m)$ of this process are independent. 
Consider
an Ito SDE in the integral form

\vspace{-1mm}
\begin{equation}
\label{1.5.2}
{\bf x}_t={\bf x}_0+\int\limits_0^t {\bf a}({\bf x}_{\tau},\tau)d\tau+
\int\limits_0^t B({\bf x}_{\tau},\tau)d{\bf f}_{\tau},\ \ \ 
{\bf x}_0={\bf x}(0,\omega).
\end{equation}

\vspace{2mm}
\noindent
Here ${\bf x}_t$ is some $n$-dimensional stochastic process 
satisfying the equation (\ref{1.5.2}). 
The nonrandom functions ${\bf a}: \mathbb{R}^n\times[0, T]\to\mathbb{R}^n$,
$B: \mathbb{R}^n\times[0, T]\to\mathbb{R}^{n\times m}$
guarantee the existence and uniqueness up to stochastic 
equivalence of a solution
of the equation (\ref{1.5.2}) \cite{1}. 
The second integral on the right-hand side of (\ref{1.5.2}) is 
interpreted as an Ito stochastic integral.
Let ${\bf x}_0$ be an $n$-dimensional random variable, which is 
${\rm F}_0$-measurable and 
${\sf M}\{\left|{\bf x}_0\right|^2\}<\infty$ 
(${\sf M}$ denotes a mathematical expectation).
We assume that
${\bf x}_0$ and ${\bf f}_t-{\bf f}_0$ are independent when $t>0.$

Consider the following iterated Ito stochastic integrals

\vspace{-1mm}
\begin{equation}
\label{ito}
J[\psi^{(k)}]_{T,t}=\int\limits_t^T\psi_k(t_k) \ldots \int\limits_t^{t_{2}}
\psi_1(t_1) d{\bf w}_{t_1}^{(i_1)}\ldots
d{\bf w}_{t_k}^{(i_k)},
\end{equation}

\vspace{2mm}
\noindent
where every $\psi_l(\tau)$ $(l=1,\ldots,k)$ is a 
non-random function 
on $[t,T],$ ${\bf w}_{\tau}^{(i)}={\bf f}_{\tau}^{(i)}$
for $i=1,\ldots,m$ and
${\bf w}_{\tau}^{(0)}=\tau,$ 
$i_1,\ldots,i_k = 0, 1,\ldots,m.$

Note that $\psi_l(\tau)\equiv 1$ $(l=1,\ldots,k)$ and
$i_1,\ldots,i_k = 0, 1,\ldots,m$ in  
\cite{KlPl2}, \cite{Mi2}, \cite{Mi3}, \cite{PW1}, \cite{KlPl1} and 
$\psi_l(\tau)\equiv (t-\tau)^{q_l}$ ($l=1,\ldots,k$; 
$q_1,\ldots,q_k=0, 1, 2,\ldots $) and $i_1,\ldots,i_k = 1,\ldots,m$ in
\cite{kk5}-\cite{2013}.
 
Effective solution 
of the problem of
combined mean-square approximation 
of the iterated Ito stochastic integrals
(\ref{ito}) of multiplicities 1 to 5
composes the subject of this article.

\vspace{5mm}

\section{Explicit One-Step Strong Numerical Schemes With Convergence
Orders 1.5, 2.0, and 2.5 Based on the Unified Taylor--Ito Expansion}

\vspace{5mm}

Introduce the definition of a strong convergence of 
a numerical method for the Ito SDEs.

Consider the partition $\{\tau_p\}_{p=0}^N$ of the interval $[0, T]$ 
with the maximum step of discretization $\Delta_N$ such that

\vspace{-4mm}
$$
0=\tau_0<\tau_1<\ldots <\tau_N=T.
$$

\vspace{2mm}

Denote by ${\bf y}_{\tau_p}\stackrel{\sf def}{=}
{\bf y}_{p},$\ $p=0, 1,\ldots,N$ the discrete approximation
of the process  ${\bf x}_t,$ $t\in[0,T]$
(the solution of the Ito SDE (\ref{1.5.2}))
that cor\-res\-ponds to the maximum discretization
step $\Delta_N$.

\vspace{2mm}

{\bf Definition 1}\ \cite{KlPl2}.\
{\it The discrete approximation {\rm (}numerical method{\rm )} 
${\bf y}_{j},$\ $j=0, 1,\ldots,N$
that cor\-res\-ponds
to the maximum discretization step $\Delta_N$
is said to be strongly 
converging with an order $\gamma>0$ 
to the process ${\bf x}_t,$ $t\in[0,T]$
if there exist a constant $C>0$ that is independent of
$\Delta_N$ and $j$ $(j=0, 1,\ldots,N)$ and also a value $\delta>0$ such that

\vspace{-1mm}
\begin{equation}
\label{yyy1}
{\sf M}\{|{\bf x}_j-{\bf y}_j|\}\le
C(\Delta_N)^{\gamma}\ \ \ (j=0, 1,\ldots,N)
\end{equation} 

\vspace{2mm}
\noindent
for all $\Delta_N\in(0, \delta).$}

\vspace{2mm}

In a series of publications \cite{Mi2}, \cite{Mi3}, the authors considered 
the mean-square convergence instead of
the strong convergence, which corresponds to the replacement 
of the con\-di\-ti\-on (\ref{yyy1}) with

\vspace{-1mm}
\begin{equation}
\label{yyy1a}
\left({\sf M}\{|{\bf x}_j-{\bf y}_j|^2\}\right)^{1/2}\le
C(\Delta_N)^{\gamma}\ \ \ (j=0, 1,\ldots,N).
\end{equation}

\vspace{2mm}

In (\ref{yyy1}) and (\ref{yyy1a}):\ ${\bf x}_{\tau_j}\stackrel{\sf def}{=}
{\bf x}_{j},$\ $j=0, 1,\ldots,N.$

Obviously, by virtue of the Lyapunov inequality \cite{KlPl2} the mean-square 
convergence implies the strong
con\-ver\-gence.
As it has appeared, a rather nontrivial question is which iterated
stochastic integrals (Ito or
Stra\-to\-no\-vich) are preferable for the numerical integration 
of the Ito SDEs with a 
correct estimation of the mean-square
approximation error. By their external view the approximations 
of the iterated Stratonovich
stochastic integrals are simpler than the corresponding 
approximations of the iterated Ito stochastic
integrals; see the details in Sections 3, 5 below. However, the estimation 
procedure of the mean-square
approximation error turns out to be much easier for the iterated
Ito stochastic integrals, which
motivates the use of these integrals.

Consider the explicit one-step numerical scheme for the 
Ito SDEs that is based on the unified
Taylor--Ito expansion \cite{2006}-\cite{2018}

\vspace{2mm}
$$
{\bf y}_{p+1}={\bf y}_p+\sum_{i_{1}=1}^{m}B_{i_{1}}
\hat I_{(0)\tau_{p+1},\tau_p}^{(i_{1})}+\Delta{\bf a}
+\sum_{i_{1},i_{2}=1}^{m}G_{i_{2}}
B_{i_{1}}\hat I_{(00)\tau_{p+1},\tau_p}^{(i_{2}i_{1})}+
$$

\vspace{2mm}
$$
+
\sum_{i_{1}=1}^{m}\Biggl(G_{i_{1}}{\bf a}\left(
\Delta \hat I_{(0)\tau_{p+1},\tau_p}^{(i_{1})}+
\hat I_{(1)\tau_{p+1},\tau_p}^{(i_{1})}\right)
-LB_{i_{1}}\hat I_{(1)\tau_{p+1},\tau_p}^{(i_{1})}\Biggr)+
$$

\vspace{2mm}
$$
+\sum_{i_{1},i_{2},i_{3}=1}^{m} G_{i_{3}}G_{i_{2}}
B_{i_{1}}\hat I_{(000)\tau_{p+1},\tau_p}^{(i_{3}i_{2}i_{1})}+
\frac{\Delta^2}{2}L{\bf a}+
$$

\vspace{1mm}
\begin{equation}
\label{4.45}
+{\bf v}_{p+1,p}+{\bf r}_{p+1,p},
\end{equation}

\vspace{4mm}
\noindent
where

\vspace{-1mm}
$$
{\bf v}_{p+1,p}=\sum_{i_{1},i_{2}=1}^{m}
\Biggl(G_{i_{2}}LB_{i_{1}}\left(
\hat I_{(10)\tau_{p+1},\tau_p}^{(i_{2}i_{1})}-
\hat I_{(01)\tau_{p+1},\tau_p}^{(i_{2}i_{1})}
\right)
-LG_{i_{2}}B_{i_{1}}\hat I_{(10)\tau_{p+1},\tau_p}^{(i_{2}i_{1})}
+\Biggr.
$$

\vspace{2mm}
$$
\Biggl.+G_{i_{2}}G_{i_{1}}{\bf a}\left(
\hat I_{(01)\tau_{p+1},\tau_p}
^{(i_{2}i_{1})}+\Delta \hat I_{(00)\tau_{p+1},\tau_p}^{(i_{2}i_{1})}
\right)\Biggr)+
$$

\vspace{2mm}
$$
+
\sum_{i_{1},i_{2},i_{3},i_{4}=1}^{m}G_{i_{4}}G_{i_{3}}G_{i_{2}}
B_{i_{1}}\hat I_{(0000)\tau_{p+1},\tau_p}^{(i_{4}i_{3}i_{2}i_{1})},
$$

\vspace{8mm}
$$
{\bf r}_{p+1,p}=\sum_{i_{1}=1}^{m}\Biggl(G_{i_{1}}L{\bf a}\left(\frac{1}{2}
\hat I_{(2)\tau_{p+1},\tau_p}
^{(i_{1})}+\Delta \hat I_{(1)\tau_{p+1},\tau_p}^{(i_{1})}+
\frac{\Delta^2}{2}\hat I_{(0)\tau_{p+1},\tau_p}^{(i_{1})}\right)\Biggr.+
$$

\vspace{2mm}
$$
\Biggl.+\frac{1}{2}LL B_{i_{1}}\hat I_{(2)\tau_{p+1},\tau_p}^{(i_{1})}-
LG_{i_{1}}{\bf a}\left(\hat I_{(2)\tau_{p+1},\tau_p}^{(i_{1})}+
\Delta \hat I_{(1)\tau_{p+1},\tau_p}^{(i_{1})}\right)\Biggr)+
$$

\vspace{2mm}
$$
+
\sum_{i_{1},i_{2},i_{3}=1}^m\Biggl(
G_{i_{3}}LG_{i_{2}}B_{i_{1}}
\left(\hat I_{(100)\tau_{p+1},\tau_p}
^{(i_{3}i_{2}i_{1})}-\hat I_{(010)\tau_{p+1},\tau_p}
^{(i_{3}i_{2}i_{1})}\right)
\Biggr.+
$$

\vspace{2mm}
$$
+G_{i_{3}}G_{i_{2}}LB_{i_{1}}\left(
\hat I_{(010)\tau_{p+1},\tau_p}^{(i_{3}i_{2}i_{1})}-
\hat I_{(001)\tau_{p+1},\tau_p}^{(i_{3}i_{2}i_{1})}\right)+
$$

\vspace{2mm}
$$
+
G_{i_{3}}G_{i_{2}}G_{i_{1}}{\bf a}
\left(\Delta \hat I_{(000)\tau_{p+1},\tau_p}^{(i_{3}i_{2}i_{1})}+
\hat I_{(001)\tau_{p+1},\tau_p}^{(i_{3}i_{2}i_{1})}\right)
-
$$

\vspace{2mm}
$$
\Biggl.-LG_{i_{3}}G_{i_{2}}B_{i_{1}}
\hat I_{(100)\tau_{p+1},\tau_p}^{(i_{3}i_{2}i_{1})}\Biggr)+
$$

\vspace{2mm}
$$
+\sum_{i_{1},i_{2},i_{3},i_{4},i_{5}=1}^m
G_{i_{5}}G_{i_{4}}G_{i_{3}}G_{i_{2}}B_{i_{1}}
\hat I_{(00000)\tau_{p+1},\tau_p}^{(i_{5}i_{4}i_{3}i_{2}i_{1})}+
$$

\vspace{2mm}
$$
+\frac{\Delta^3}{6}LL{\bf a},
$$

\vspace{6mm}
\noindent
where $\Delta=T/N$ $(N>1)$ is a constant (for simplicity) 
integration step,\
$\tau_p=p\Delta$ $(p=0, 1,\ldots,N)$,\
$\hat I_{(l_1\ldots\hspace{0.2mm} l_k)s,t}^{(i_1\ldots i_k)}$ 
denotes
an approximation of the iterated Ito stochastic integral of multiplicity
$k$

\vspace{-1mm}
\begin{equation}
\label{ll1}
I_{(l_1\ldots\hspace{0.2mm} l_k)s,t}^{(i_1\ldots i_k)}=
 \int\limits^ {s} _ {t} (t-\tau _
{k}) ^ {l_ {k}} 
\ldots \int\limits^ {\tau _ {2}} _ {t} (t-\tau _ {1}) ^ {l_ {1}} d
{\bf f} ^ {(i_ {1})} _ {\tau_ {1}} \ldots 
d {\bf f} _ {\tau_ {k}} ^ {(i_ {k})},
\end{equation}

\vspace{3mm}
$$
L= {\partial \over \partial t}
+ \sum^ {n} _ {i=1} {\bf a}_i ({\bf x},  t) 
{\partial  \over  \partial  {\bf  x}_i}
+ {1\over 2} \sum^ {m} _ {j=1} \sum^ {n} _ {l,i=1}
B_{lj} ({\bf x}, t) B_{ij} ({\bf x}, t) {\partial
^{2} \over \partial {\bf x}_l \partial {\bf x}_i},
$$

\vspace{3mm}
$$
G_i = \sum^ {n} _ {j=1} B_{ji} ({\bf x}, t)
{\partial  \over \partial {\bf x}_j},\ \ \ 
i=1,\ldots,m,
$$

\vspace{3mm}
\noindent
$l_1,\ldots, l_k=0, 1, 2,$\ \
$i_1,\ldots, i_k=1,\ldots,m,$\ \ $k=1, 2,\ldots, 5$,\ \
$B_i$ and $B_{ij}$ are the $i$th column and the $ij$th
element of the matrix function $B$,
${\bf a}_i$ is the $i$th element of the vector function ${\bf a},$
${\bf x}_i$ is the $i$th element
of the column ${\bf x}$, 
the functions  

\vspace{-1mm}
$$
B_{i_{1}},\ {\bf a},\ G_{i_{2}}B_{i_{1}},\
G_{i_{1}}{\bf a},\ LB_{i_{1}},\ G_{i_{3}}G_{i_{2}}B_{i_{1}},\ 
L{\bf a},\ LL{\bf a},\
G_{i_{2}}LB_{i_{1}},
$$

\vspace{-3mm}
$$
LG_{i_{2}}B_{i_{1}},\ G_{i_{2}}G_{i_{1}}{\bf a},\
G_{i_{4}}G_{i_{3}}G_{i_{2}}B_{i_{1}},\ G_{i_{1}}L{\bf a},\
LLB_{i_{1}},\ LG_{i_{1}}{\bf a},\ G_{i_{3}}LG_{i_{2}}B_{i_{1}},\
G_{i_{3}}G_{i_{2}}LB_{i_{1}},
$$

\vspace{-3mm}
$$
G_{i_{3}}G_{i_{2}}G_{i_{1}}{\bf a},\
LG_{i_{3}}G_{i_{2}}B_{i_{1}},\ 
G_{i_{5}}G_{i_{4}}G_{i_{3}}G_{i_{2}}B_{i_{1}}
$$

\vspace{3mm}
\noindent
are calculated at the point $({\bf y}_p,p).$

Under the standard conditions \cite{KlPl2}, \cite{2006} the numerical 
scheme (\ref{4.45}) has 
order 2.5 of strong convergence. The major emphasis below will be placed on the 
approximation of the iterated
Ito stochastic integrals appearing in (\ref{4.45}). Therefore, among 
the standard conditions, note the
approximation condi\-ti\-on of these integrals \cite{KlPl2}, \cite{2006}, 
which has the form

\vspace{-1mm}
\begin{equation}
\label{4.3}
{\sf M}\left\{\biggl(I_{(l_{1}\ldots\hspace{0.2mm} l_{k})\tau_{p+1},\tau_p}
^{(i_{1}\ldots i_{k})} 
-\hat I_{(l_{1}
\ldots\hspace{0.2mm} l_{k})\tau_{p+1},\tau_p}^{(i_{1}\ldots i_{k})}
\biggr)^2\right\}\le C\Delta^{6},
\end{equation}

\vspace{2mm}
\noindent
where constant $C$ is independent of
$\Delta$.

Note that if we exclude ${\bf v}_{p+1,p}+{\bf r}_{p+1,p}$ from the
right-hand side of (\ref{4.45}), then we have an explicit 
one-step strong numerical scheme of order 1.5 \cite{KlPl2}, 
\cite{2006}-\cite{2018}.
The right-hand side of (\ref{4.45}) but without the value
${\bf r}_{p+1,p}$ define 
an explicit 
one-step strong numerical scheme of order 2.0
\cite{KlPl2}, 
\cite{2006}-\cite{2018}.

Using the numerical scheme (\ref{4.45}) or its modification based 
on the Taylor--Ito expansion \cite{KlPl1},
the implicit or multistep analogs of (\ref{4.45}) can be constructed; 
see \cite{KlPl2}. The set of the
iterated Ito stochastic integrals to be approximated for implementing 
these modifications is the same
as for the numerical scheme (\ref{4.45}) itself.
Interestingly, the truncated unified Taylor--Ito expansion --- the 
foundation of the numerical
scheme (\ref{4.45}) --- contains 12 different iterated Ito 
stochastic integrals 
of the form (\ref{ll1}), which cannot be
interconnected by linear relations \cite{2006}-\cite{2018}. The analogous 
Taylor--Ito expansion \cite{KlPl2} contains
17 different iterated Ito stochastic integrals, part of which 
are interconnected by linear relations
and part of which have a higher multiplicity than the iterated 
Ito stochastic integrals (\ref{ll1}). This
fact well explains the use of the numerical scheme (\ref{4.45}).

One of the main problems arising in the implementation of the 
numerical scheme (\ref{4.45}) is the joint
numerical modeling of the iterated Ito stochastic integrals 
figuring in (\ref{4.45}). In the next section,
we will consider an efficient numerical modeling method for 
the iterated Ito stochastic integrals
and also demonstrate which stochastic integrals (Ito or Stratonovich) 
are preferable for numerical
modeling with a correct estimation of the mean-square approximation error.

\vspace{5mm}

\section{Method of Numerical Modeling for Iterated Ito 
Stochastic Integrals
Based on Genegalized Multiple Fourier Series}

\vspace{5mm}

An efficient numerical modeling method for the iterated Ito 
stochastic integrals based on generalized multiple
Fourier series was considered in \cite{2006} (2006); also see 
\cite{2011-2}-\cite{new-art-1-xxy}. 
This method rests on an important
result presented below (Theorems 1, 2).

Suppose that $\psi_1(\tau),\ldots,\psi_k(\tau)\in L_2([t, T])$.
Define the following function on the hypercube $[t, T]^k$

\vspace{-1mm}
\begin{equation}
\label{ppp}
K(t_1,\ldots,t_k)=
\begin{cases}
\psi_1(t_1)\ldots \psi_k(t_k)\ &\hbox{for}\ \ t_1<\ldots<t_k\\
~\\
~\\
0\ &\hbox{otherwise}
\end{cases},\ \ \ \ t_1,\ldots,t_k\in[t, T],\ \ \ \ k\ge 2,
\end{equation}

\vspace{4mm}
\noindent
and 
$K(t_1)\equiv\psi_1(t_1)$ for $t_1\in[t, T].$

Suppose that $\{\phi_j(x)\}_{j=0}^{\infty}$
is a complete orthonormal system of functions in the space
$L_2([t, T])$. 
The function $K(t_1,\ldots,t_k)$ belongs to the space $L_2([t, T]^k).$
At this situation it is well known that the generalized 
multiple Fourier series 
of $K(t_1,\ldots,t_k)\in L_2([t, T]^k)$ is converging 
to $K(t_1,\ldots,t_k)$ in the hypercube $[t, T]^k$ in 
the mean-square sense, i.e.

$$
\hbox{\vtop{\offinterlineskip\halign{
\hfil#\hfil\cr
{\rm lim}\cr
$\stackrel{}{{}_{p_1,\ldots,p_k\to \infty}}$\cr
}} }\left\Vert
K(t_1,\ldots,t_k)-
\sum_{j_1=0}^{p_1}\ldots \sum_{j_k=0}^{p_k}
C_{j_k\ldots j_1}\prod_{l=1}^{k} \phi_{j_l}(t_l)
\right\Vert_{L_2([t,T]^k)}=0,
$$

\vspace{3mm}
\noindent
where
\begin{equation}
\label{ppppa}
C_{j_k\ldots j_1}=\int\limits_{[t,T]^k}
K(t_1,\ldots,t_k)\prod_{l=1}^{k}\phi_{j_l}(t_l)dt_1\ldots dt_k,
\end{equation}

$$
\left\Vert f\right\Vert_{L_2([t,T]^k)}=\left(\int\limits_{[t,T]^k}
f^2(t_1,\ldots,t_k)dt_1\ldots dt_k\right)^{1/2}.
$$

\vspace{5mm}

Consider the partition $\{\tau_j\}_{j=0}^N$ of $[t,T]$ such that

\vspace{-1mm}
\begin{equation}
\label{1111}
t=\tau_0<\ldots <\tau_N=T,\ \ \
\Delta_N=
\hbox{\vtop{\offinterlineskip\halign{
\hfil#\hfil\cr
{\rm max}\cr
$\stackrel{}{{}_{0\le j\le N-1}}$\cr
}} }\Delta\tau_j\to 0\ \ \hbox{if}\ \ N\to \infty,\ \ \
\Delta\tau_j=\tau_{j+1}-\tau_j.
\end{equation}

\vspace{3mm}                

{\bf Theorem 1}\ \cite{2006} (2006),
\cite{2011-2}-\cite{200a}, \cite{301a}-\cite{arxiv}, 
\cite{arxiv-24}-\cite{new-art-1-xxy}.
{\it Suppose that
every 
$\psi_l(\tau)$ $(l=1,\ldots, k)$ is a continuous non-random func\-tion on 
$[t, T]$ and
$\{\phi_j(x)\}_{j=0}^{\infty}$ is a complete orthonormal system  
of continuous func\-ti\-ons in the space $L_2([t,T]).$ Then

$$
J[\psi^{(k)}]_{T,t}\  =\ 
\hbox{\vtop{\offinterlineskip\halign{
\hfil#\hfil\cr
{\rm l.i.m.}\cr
$\stackrel{}{{}_{p_1,\ldots,p_k\to \infty}}$\cr
}} }\sum_{j_1=0}^{p_1}\ldots\sum_{j_k=0}^{p_k}
C_{j_k\ldots j_1}\Biggl(
\prod_{l=1}^k\zeta_{j_l}^{(i_l)}\ -
\Biggr.
$$

\vspace{2mm}
\begin{equation}
\label{tyyy}
-\ \Biggl.
\hbox{\vtop{\offinterlineskip\halign{
\hfil#\hfil\cr
{\rm l.i.m.}\cr
$\stackrel{}{{}_{N\to \infty}}$\cr
}} }\sum_{(l_1,\ldots,l_k)\in {\rm G}_k}
\phi_{j_{1}}(\tau_{l_1})
\Delta{\bf w}_{\tau_{l_1}}^{(i_1)}\ldots
\phi_{j_{k}}(\tau_{l_k})
\Delta{\bf w}_{\tau_{l_k}}^{(i_k)}\Biggr),
\end{equation}

\vspace{5mm}
\noindent
where $J[\psi^{(k)}]_{T,t}$ is defined by {\rm (\ref{ito}),}

\vspace{-1mm}
$$
{\rm G}_k={\rm H}_k\backslash{\rm L}_k,\ \ \
{\rm H}_k=\{(l_1,\ldots,l_k):\ l_1,\ldots,l_k=0,\ 1,\ldots,N-1\},
$$

\vspace{-1mm}
$$
{\rm L}_k=\{(l_1,\ldots,l_k):\ l_1,\ldots,l_k=0,\ 1,\ldots,N-1;\
l_g\ne l_r\ (g\ne r);\ g, r=1,\ldots,k\},
$$

\vspace{3mm}
\noindent
${\rm l.i.m.}$ is a limit in the mean-square sense$,$
$i_1,\ldots,i_k=0,1,\ldots,m,$

\vspace{-1mm}
\begin{equation}
\label{rr23}
\zeta_{j}^{(i)}=
\int\limits_t^T \phi_{j}(s) d{\bf w}_s^{(i)}
\end{equation} 

\vspace{2mm}
\noindent
are independent standard Gaussian random variables
for various
$i$ or $j$ {\rm(}if $i\ne 0${\rm),}
$C_{j_k\ldots j_1}$ is the Fourier coefficient {\rm(\ref{ppppa}),}
$\Delta{\bf w}_{\tau_{j}}^{(i)}=
{\bf w}_{\tau_{j+1}}^{(i)}-{\bf w}_{\tau_{j}}^{(i)}$
$(i=0, 1,\ldots,m),$
$\left\{\tau_{j}\right\}_{j=0}^{N}$ is a partition of
the interval $[t, T],$ which satisfies the condition {\rm (\ref{1111})}.
}

\vspace{2mm}

The convergence in the mean of degree 
$2n$ $(n\in \mathbb{N})$ \cite{2018a}-\cite{2018aa1}, 
\cite{2007-2}-\cite{2013}
as well as the convergence
with probability 1 \cite{2018a}-\cite{2018aa1}, \cite{arxiv-1},
\cite{arxiv-3}
of approximations from Theorem 1 (also see Theorem 2 below)
are proved.
Moreover, the complete orthonormal systems of Haar and 
Rademacher--Walsh functions in $L_2([t,T])$ 
can also be applied in Theorems 1
\cite{2006}-\cite{2013}.
The modification of Theorem 1 for 
complete orthonormal with weigth $r(x)\ge 0$ systems
of functions in the space $L_2([t,T])$ can be found in 
\cite{2018}-\cite{2018aa1}, \cite{arxiv-26b}

In order to evaluate the significance of Theorem 1 for practice we will
demonstrate its transfor\-med particular cases for 
$k=1,\ldots,5$ 
\cite{2006}-\cite{200a}, \cite{301a}-\cite{arxiv}, 
\cite{arxiv-24}-\cite{new-art-1-xxy}

\vspace{1mm}
\begin{equation}
\label{a1}
J[\psi^{(1)}]_{T,t}
=\hbox{\vtop{\offinterlineskip\halign{
\hfil#\hfil\cr
{\rm l.i.m.}\cr
$\stackrel{}{{}_{p_1\to \infty}}$\cr
}} }\sum_{j_1=0}^{p_1}
C_{j_1}\zeta_{j_1}^{(i_1)},
\end{equation}

\vspace{2mm}
\begin{equation}
\label{a2}
J[\psi^{(2)}]_{T,t}
=\hbox{\vtop{\offinterlineskip\halign{
\hfil#\hfil\cr
{\rm l.i.m.}\cr
$\stackrel{}{{}_{p_1,p_2\to \infty}}$\cr
}} }\sum_{j_1=0}^{p_1}\sum_{j_2=0}^{p_2}
C_{j_2j_1}\Biggl(\zeta_{j_1}^{(i_1)}\zeta_{j_2}^{(i_2)}
-{\bf 1}_{\{i_1=i_2\ne 0\}}
{\bf 1}_{\{j_1=j_2\}}\Biggr),
\end{equation}

\vspace{6mm}

$$
J[\psi^{(3)}]_{T,t}=
\hbox{\vtop{\offinterlineskip\halign{
\hfil#\hfil\cr
{\rm l.i.m.}\cr
$\stackrel{}{{}_{p_1,p_2,p_3\to \infty}}$\cr
}} }\sum_{j_1=0}^{p_1}\sum_{j_2=0}^{p_2}\sum_{j_3=0}^{p_3}
C_{j_3j_2j_1}\Biggl(
\zeta_{j_1}^{(i_1)}\zeta_{j_2}^{(i_2)}\zeta_{j_3}^{(i_3)}
-\Biggr.
$$

\begin{equation}
\label{a3}
-\Biggl.
{\bf 1}_{\{i_1=i_2\ne 0\}}
{\bf 1}_{\{j_1=j_2\}}
\zeta_{j_3}^{(i_3)}
-{\bf 1}_{\{i_2=i_3\ne 0\}}
{\bf 1}_{\{j_2=j_3\}}
\zeta_{j_1}^{(i_1)}-
{\bf 1}_{\{i_1=i_3\ne 0\}}
{\bf 1}_{\{j_1=j_3\}}
\zeta_{j_2}^{(i_2)}\Biggr),
\end{equation}

\vspace{6mm}

$$
J[\psi^{(4)}]_{T,t}
=
\hbox{\vtop{\offinterlineskip\halign{
\hfil#\hfil\cr
{\rm l.i.m.}\cr
$\stackrel{}{{}_{p_1,\ldots,p_4\to \infty}}$\cr
}} }\sum_{j_1=0}^{p_1}\ldots\sum_{j_4=0}^{p_4}
C_{j_4\ldots j_1}\Biggl(
\prod_{l=1}^4\zeta_{j_l}^{(i_l)}
\Biggr.
-
$$

\vspace{-1mm}
$$
-
{\bf 1}_{\{i_1=i_2\ne 0\}}
{\bf 1}_{\{j_1=j_2\}}
\zeta_{j_3}^{(i_3)}
\zeta_{j_4}^{(i_4)}
-
{\bf 1}_{\{i_1=i_3\ne 0\}}
{\bf 1}_{\{j_1=j_3\}}
\zeta_{j_2}^{(i_2)}
\zeta_{j_4}^{(i_4)}-
$$

\vspace{-2.5mm}
$$
-
{\bf 1}_{\{i_1=i_4\ne 0\}}
{\bf 1}_{\{j_1=j_4\}}
\zeta_{j_2}^{(i_2)}
\zeta_{j_3}^{(i_3)}
-
{\bf 1}_{\{i_2=i_3\ne 0\}}
{\bf 1}_{\{j_2=j_3\}}
\zeta_{j_1}^{(i_1)}
\zeta_{j_4}^{(i_4)}-
$$

\vspace{-2.5mm}
$$
-
{\bf 1}_{\{i_2=i_4\ne 0\}}
{\bf 1}_{\{j_2=j_4\}}
\zeta_{j_1}^{(i_1)}
\zeta_{j_3}^{(i_3)}
-
{\bf 1}_{\{i_3=i_4\ne 0\}}
{\bf 1}_{\{j_3=j_4\}}
\zeta_{j_1}^{(i_1)}
\zeta_{j_2}^{(i_2)}+
$$

\vspace{-2.5mm}
$$
+
{\bf 1}_{\{i_1=i_2\ne 0\}}
{\bf 1}_{\{j_1=j_2\}}
{\bf 1}_{\{i_3=i_4\ne 0\}}
{\bf 1}_{\{j_3=j_4\}}
+
$$

\vspace{-2.5mm}
$$
+
{\bf 1}_{\{i_1=i_3\ne 0\}}
{\bf 1}_{\{j_1=j_3\}}
{\bf 1}_{\{i_2=i_4\ne 0\}}
{\bf 1}_{\{j_2=j_4\}}+
$$
\begin{equation}
\label{a4}
+\Biggl.
{\bf 1}_{\{i_1=i_4\ne 0\}}
{\bf 1}_{\{j_1=j_4\}}
{\bf 1}_{\{i_2=i_3\ne 0\}}
{\bf 1}_{\{j_2=j_3\}}\Biggr),
\end{equation}

\vspace{6mm}

$$
J[\psi^{(5)}]_{T,t}
=\hbox{\vtop{\offinterlineskip\halign{
\hfil#\hfil\cr
{\rm l.i.m.}\cr
$\stackrel{}{{}_{p_1,\ldots,p_5\to \infty}}$\cr
}} }\sum_{j_1=0}^{p_1}\ldots\sum_{j_5=0}^{p_5}
C_{j_5\ldots j_1}\Biggl(
\prod_{l=1}^5\zeta_{j_l}^{(i_l)}
-\Biggr.
$$
$$
-
{\bf 1}_{\{i_1=i_2\ne 0\}}
{\bf 1}_{\{j_1=j_2\}}
\zeta_{j_3}^{(i_3)}
\zeta_{j_4}^{(i_4)}
\zeta_{j_5}^{(i_5)}-
{\bf 1}_{\{i_1=i_3\ne 0\}}
{\bf 1}_{\{j_1=j_3\}}
\zeta_{j_2}^{(i_2)}
\zeta_{j_4}^{(i_4)}
\zeta_{j_5}^{(i_5)}-
$$
$$
-
{\bf 1}_{\{i_1=i_4\ne 0\}}
{\bf 1}_{\{j_1=j_4\}}
\zeta_{j_2}^{(i_2)}
\zeta_{j_3}^{(i_3)}
\zeta_{j_5}^{(i_5)}-
{\bf 1}_{\{i_1=i_5\ne 0\}}
{\bf 1}_{\{j_1=j_5\}}
\zeta_{j_2}^{(i_2)}
\zeta_{j_3}^{(i_3)}
\zeta_{j_4}^{(i_4)}-
$$
$$
-
{\bf 1}_{\{i_2=i_3\ne 0\}}
{\bf 1}_{\{j_2=j_3\}}
\zeta_{j_1}^{(i_1)}
\zeta_{j_4}^{(i_4)}
\zeta_{j_5}^{(i_5)}-
{\bf 1}_{\{i_2=i_4\ne 0\}}
{\bf 1}_{\{j_2=j_4\}}
\zeta_{j_1}^{(i_1)}
\zeta_{j_3}^{(i_3)}
\zeta_{j_5}^{(i_5)}-
$$
$$
-
{\bf 1}_{\{i_2=i_5\ne 0\}}
{\bf 1}_{\{j_2=j_5\}}
\zeta_{j_1}^{(i_1)}
\zeta_{j_3}^{(i_3)}
\zeta_{j_4}^{(i_4)}
-{\bf 1}_{\{i_3=i_4\ne 0\}}
{\bf 1}_{\{j_3=j_4\}}
\zeta_{j_1}^{(i_1)}
\zeta_{j_2}^{(i_2)}
\zeta_{j_5}^{(i_5)}-
$$
$$
-
{\bf 1}_{\{i_3=i_5\ne 0\}}
{\bf 1}_{\{j_3=j_5\}}
\zeta_{j_1}^{(i_1)}
\zeta_{j_2}^{(i_2)}
\zeta_{j_4}^{(i_4)}
-{\bf 1}_{\{i_4=i_5\ne 0\}}
{\bf 1}_{\{j_4=j_5\}}
\zeta_{j_1}^{(i_1)}
\zeta_{j_2}^{(i_2)}
\zeta_{j_3}^{(i_3)}+
$$
$$
+
{\bf 1}_{\{i_1=i_2\ne 0\}}
{\bf 1}_{\{j_1=j_2\}}
{\bf 1}_{\{i_3=i_4\ne 0\}}
{\bf 1}_{\{j_3=j_4\}}\zeta_{j_5}^{(i_5)}+
{\bf 1}_{\{i_1=i_2\ne 0\}}
{\bf 1}_{\{j_1=j_2\}}
{\bf 1}_{\{i_3=i_5\ne 0\}}
{\bf 1}_{\{j_3=j_5\}}\zeta_{j_4}^{(i_4)}+
$$
$$
+
{\bf 1}_{\{i_1=i_2\ne 0\}}
{\bf 1}_{\{j_1=j_2\}}
{\bf 1}_{\{i_4=i_5\ne 0\}}
{\bf 1}_{\{j_4=j_5\}}\zeta_{j_3}^{(i_3)}+
{\bf 1}_{\{i_1=i_3\ne 0\}}
{\bf 1}_{\{j_1=j_3\}}
{\bf 1}_{\{i_2=i_4\ne 0\}}
{\bf 1}_{\{j_2=j_4\}}\zeta_{j_5}^{(i_5)}+
$$
$$
+
{\bf 1}_{\{i_1=i_3\ne 0\}}
{\bf 1}_{\{j_1=j_3\}}
{\bf 1}_{\{i_2=i_5\ne 0\}}
{\bf 1}_{\{j_2=j_5\}}\zeta_{j_4}^{(i_4)}+
{\bf 1}_{\{i_1=i_3\ne 0\}}
{\bf 1}_{\{j_1=j_3\}}
{\bf 1}_{\{i_4=i_5\ne 0\}}
{\bf 1}_{\{j_4=j_5\}}\zeta_{j_2}^{(i_2)}+
$$
$$
+
{\bf 1}_{\{i_1=i_4\ne 0\}}
{\bf 1}_{\{j_1=j_4\}}
{\bf 1}_{\{i_2=i_3\ne 0\}}
{\bf 1}_{\{j_2=j_3\}}\zeta_{j_5}^{(i_5)}+
{\bf 1}_{\{i_1=i_4\ne 0\}}
{\bf 1}_{\{j_1=j_4\}}
{\bf 1}_{\{i_2=i_5\ne 0\}}
{\bf 1}_{\{j_2=j_5\}}\zeta_{j_3}^{(i_3)}+
$$
$$
+
{\bf 1}_{\{i_1=i_4\ne 0\}}
{\bf 1}_{\{j_1=j_4\}}
{\bf 1}_{\{i_3=i_5\ne 0\}}
{\bf 1}_{\{j_3=j_5\}}\zeta_{j_2}^{(i_2)}+
{\bf 1}_{\{i_1=i_5\ne 0\}}
{\bf 1}_{\{j_1=j_5\}}
{\bf 1}_{\{i_2=i_3\ne 0\}}
{\bf 1}_{\{j_2=j_3\}}\zeta_{j_4}^{(i_4)}+
$$
$$
+
{\bf 1}_{\{i_1=i_5\ne 0\}}
{\bf 1}_{\{j_1=j_5\}}
{\bf 1}_{\{i_2=i_4\ne 0\}}
{\bf 1}_{\{j_2=j_4\}}\zeta_{j_3}^{(i_3)}+
{\bf 1}_{\{i_1=i_5\ne 0\}}
{\bf 1}_{\{j_1=j_5\}}
{\bf 1}_{\{i_3=i_4\ne 0\}}
{\bf 1}_{\{j_3=j_4\}}\zeta_{j_2}^{(i_2)}+
$$
$$
+
{\bf 1}_{\{i_2=i_3\ne 0\}}
{\bf 1}_{\{j_2=j_3\}}
{\bf 1}_{\{i_4=i_5\ne 0\}}
{\bf 1}_{\{j_4=j_5\}}\zeta_{j_1}^{(i_1)}+
{\bf 1}_{\{i_2=i_4\ne 0\}}
{\bf 1}_{\{j_2=j_4\}}
{\bf 1}_{\{i_3=i_5\ne 0\}}
{\bf 1}_{\{j_3=j_5\}}\zeta_{j_1}^{(i_1)}+
$$
\begin{equation}
\label{a5}
+\Biggl.
{\bf 1}_{\{i_2=i_5\ne 0\}}
{\bf 1}_{\{j_2=j_5\}}
{\bf 1}_{\{i_3=i_4\ne 0\}}
{\bf 1}_{\{j_3=j_4\}}\zeta_{j_1}^{(i_1)}\Biggr),
\end{equation}

\vspace{5mm}
\noindent
where ${\bf 1}_A$ is the indicator of the set $A$.

We will consider the case $i_1,\ldots,i_5=1,\ldots,m$.
This case corresponds to the numerical method (\ref{4.45}).

For further consideration, let us 
consider the generalization of formulas (\ref{a1})--(\ref{a5})                 
for the case of an arbitrary multiplicity $k$ $(k\in\mathbb{N})$ of 
the iterated Ito stochastic integral $J[\psi^{(k)}]_{T,t}$ defined by (\ref{ito}).
In order to do this, let us
introduce some notations. 
Consider the unordered
set $\{1, 2, \ldots, k\}$ 
and separate it into two parts:
the first part consists of $r$ unordered 
pairs (sequence order of these pairs is also unimportant) and the 
second one consists of the 
remaining $k-2r$ numbers.
So, we have

\begin{equation}
\label{leto5007}
(\{
\underbrace{\{g_1, g_2\}, \ldots, 
\{g_{2r-1}, g_{2r}\}}_{\small{\hbox{part 1}}}
\},
\{\underbrace{q_1, \ldots, q_{k-2r}}_{\small{\hbox{part 2}}}
\}),
\end{equation}

\vspace{4mm}
\noindent
where 

\vspace{-2mm}
$$
\{g_1, g_2, \ldots, 
g_{2r-1}, g_{2r}, q_1, \ldots, q_{k-2r}\}=\{1, 2, \ldots, k\},
$$

\vspace{4mm}
\noindent
braces   
mean an unordered 
set, and pa\-ren\-the\-ses mean an ordered set.

We will say that (\ref{leto5007}) is a partition 
and consider the sum with respect to all possible
partitions

\begin{equation}
\label{leto5008}
\sum_{\stackrel{(\{\{g_1, g_2\}, \ldots, 
\{g_{2r-1}, g_{2r}\}\}, \{q_1, \ldots, q_{k-2r}\})}
{{}_{\{g_1, g_2, \ldots, 
g_{2r-1}, g_{2r}, q_1, \ldots, q_{k-2r}\}=\{1, 2, \ldots, k\}}}}
a_{g_1 g_2, \ldots, 
g_{2r-1} g_{2r}, q_1 \ldots q_{k-2r}}.
\end{equation}

\vspace{4mm}

Below there are several examples of sums in the form (\ref{leto5008})

\vspace{2mm}
$$
\sum_{\stackrel{(\{g_1, g_2\})}{{}_{\{g_1, g_2\}=\{1, 2\}}}}
a_{g_1 g_2}=a_{12},
$$

\vspace{3mm}
$$
\sum_{\stackrel{(\{\{g_1, g_2\}, \{g_3, g_4\}\})}
{{}_{\{g_1, g_2, g_3, g_4\}=\{1, 2, 3, 4\}}}}
a_{g_1 g_2 g_3 g_4}=a_{1234} + a_{1324} + a_{2314},
$$

\vspace{3mm}
$$
\sum_{\stackrel{(\{g_1, g_2\}, \{q_1, q_{2}\})}
{{}_{\{g_1, g_2, q_1, q_{2}\}=\{1, 2, 3, 4\}}}}
a_{g_1 g_2, q_1 q_{2}}=
$$

$$
=a_{12,34}+a_{13,24}+a_{14,23}
+a_{23,14}+a_{24,13}+a_{34,12},
$$

\vspace{3mm}
$$
\sum_{\stackrel{(\{g_1, g_2\}, \{q_1, q_{2}, q_3\})}
{{}_{\{g_1, g_2, q_1, q_{2}, q_3\}=\{1, 2, 3, 4, 5\}}}}
a_{g_1 g_2, q_1 q_{2}q_3}
=
$$

$$
=a_{12,345}+a_{13,245}+a_{14,235}
+a_{15,234}+a_{23,145}+a_{24,135}+
$$
$$
+a_{25,134}+a_{34,125}+a_{35,124}+a_{45,123},
$$

\vspace{4mm}
$$
\sum_{\stackrel{(\{\{g_1, g_2\}, \{g_3, g_{4}\}\}, \{q_1\})}
{{}_{\{g_1, g_2, g_3, g_{4}, q_1\}=\{1, 2, 3, 4, 5\}}}}
a_{g_1 g_2, g_3 g_{4},q_1}
=
$$

$$
=
a_{12,34,5}+a_{13,24,5}+a_{14,23,5}+
a_{12,35,4}+a_{13,25,4}+a_{15,23,4}+
$$
$$
+a_{12,54,3}+a_{15,24,3}+a_{14,25,3}+a_{15,34,2}+a_{13,54,2}+a_{14,53,2}+
$$
$$
+
a_{52,34,1}+a_{53,24,1}+a_{54,23,1}.
$$

\vspace{5mm}

Now we can write (\ref{tyyy}) as

\vspace{1mm}

$$
J[\psi^{(k)}]_{T,t}=
\hbox{\vtop{\offinterlineskip\halign{
\hfil#\hfil\cr
{\rm l.i.m.}\cr
$\stackrel{}{{}_{p_1,\ldots,p_k\to \infty}}$\cr
}} }
\sum\limits_{j_1=0}^{p_1}\ldots
\sum\limits_{j_k=0}^{p_k}
C_{j_k\ldots j_1}\Biggl(
\prod_{l=1}^k\zeta_{j_l}^{(i_l)}+\sum\limits_{r=1}^{[k/2]}
(-1)^r \times
\Biggr.
$$

\vspace{3mm}
\begin{equation}
\label{leto6000hh}
\times
\sum_{\stackrel{(\{\{g_1, g_2\}, \ldots, 
\{g_{2r-1}, g_{2r}\}\}, \{q_1, \ldots, q_{k-2r}\})}
{{}_{\{g_1, g_2, \ldots, 
g_{2r-1}, g_{2r}, q_1, \ldots, q_{k-2r}\}=\{1, 2, \ldots, k\}}}}
\prod\limits_{s=1}^r
{\bf 1}_{\{i_{g_{{}_{2s-1}}}=~i_{g_{{}_{2s}}}\ne 0\}}
\Biggl.{\bf 1}_{\{j_{g_{{}_{2s-1}}}=~j_{g_{{}_{2s}}}\}}
\prod_{l=1}^{k-2r}\zeta_{j_{q_l}}^{(i_{q_l})}\Biggr),
\end{equation}

\vspace{5mm}
\noindent
where $[x]$ is an integer part of a real number $x;$
another notations are the same as in Theorem {\bf 1}.

\vspace{2mm}

In particular, from (\ref{leto6000hh}) for $k=5$ we obtain

\vspace{3mm}

$$
J[\psi^{(5)}]_{T,t}=
\hbox{\vtop{\offinterlineskip\halign{
\hfil#\hfil\cr
{\rm l.i.m.}\cr
$\stackrel{}{{}_{p_1,\ldots,p_5\to \infty}}$\cr
}} }\sum_{j_1=0}^{p_1}\ldots\sum_{j_5=0}^{p_5}
C_{j_5\ldots j_1}\Biggl(
\prod_{l=1}^5\zeta_{j_l}^{(i_l)}-\Biggr.
$$

\vspace{2mm}
$$
-
\sum\limits_{\stackrel{(\{g_1, g_2\}, \{q_1, q_{2}, q_3\})}
{{}_{\{g_1, g_2, q_{1}, q_{2}, q_3\}=\{1, 2, 3, 4, 5\}}}}
{\bf 1}_{\{i_{g_{{}_{1}}}=~i_{g_{{}_{2}}}\ne 0\}}
{\bf 1}_{\{j_{g_{{}_{1}}}=~j_{g_{{}_{2}}}\}}
\prod_{l=1}^{3}\zeta_{j_{q_l}}^{(i_{q_l})}+
$$

\vspace{2mm}
$$
+
\sum_{\stackrel{(\{\{g_1, g_2\}, 
\{g_{3}, g_{4}\}\}, \{q_1\})}
{{}_{\{g_1, g_2, g_{3}, g_{4}, q_1\}=\{1, 2, 3, 4, 5\}}}}
{\bf 1}_{\{i_{g_{{}_{1}}}=~i_{g_{{}_{2}}}\ne 0\}}
{\bf 1}_{\{j_{g_{{}_{1}}}=~j_{g_{{}_{2}}}\}}
\Biggl.{\bf 1}_{\{i_{g_{{}_{3}}}=~i_{g_{{}_{4}}}\ne 0\}}
{\bf 1}_{\{j_{g_{{}_{3}}}=~j_{g_{{}_{4}}}\}}
\zeta_{j_{q_1}}^{(i_{q_1})}\Biggr).
$$

\vspace{7mm}
\noindent
The last equality obviously agrees with
(\ref{a5}).

Let us consider the generalization of Theorem 1 for the case
of an arbitrary complete orthonormal systems  
of functions in the space $L_2([t,T])$ 
and $\psi_1(\tau),\ldots,\psi_k(\tau)\in L_2([t, T]).$

\vspace{2mm}

{\bf Theorem~2}\ \cite{2018a} (Sect.~1.11), \cite{arxiv-1} (Sect.~15).
{\it Suppose that
$\psi_1(\tau),\ldots,\psi_k(\tau)\in L_2([t, T])$ and
$\{\phi_j(x)\}_{j=0}^{\infty}$ is an arbitrary complete orthonormal system  
of functions in the space $L_2([t,T]).$
Then the following expansion

\vspace{1mm}
$$
J[\psi^{(k)}]_{T,t}=
\hbox{\vtop{\offinterlineskip\halign{
\hfil#\hfil\cr
{\rm l.i.m.}\cr
$\stackrel{}{{}_{p_1,\ldots,p_k\to \infty}}$\cr
}} }
\sum\limits_{j_1=0}^{p_1}\ldots
\sum\limits_{j_k=0}^{p_k}
C_{j_k\ldots j_1}\Biggl(
\prod_{l=1}^k\zeta_{j_l}^{(i_l)}+\sum\limits_{r=1}^{[k/2]}
(-1)^r \times
\Biggr.
$$

\vspace{2mm}
\begin{equation}
\label{leto6000}
\times
\sum_{\stackrel{(\{\{g_1, g_2\}, \ldots, 
\{g_{2r-1}, g_{2r}\}\}, \{q_1, \ldots, q_{k-2r}\})}
{{}_{\{g_1, g_2, \ldots, 
g_{2r-1}, g_{2r}, q_1, \ldots, q_{k-2r}\}=\{1, 2, \ldots, k\}}}}
\prod\limits_{s=1}^r
{\bf 1}_{\{i_{g_{{}_{2s-1}}}=~i_{g_{{}_{2s}}}\ne 0\}}
\Biggl.{\bf 1}_{\{j_{g_{{}_{2s-1}}}=~j_{g_{{}_{2s}}}\}}
\prod_{l=1}^{k-2r}\zeta_{j_{q_l}}^{(i_{q_l})}\Biggr)
\end{equation}

\vspace{6mm}
\noindent
con\-verg\-ing in the mean-square sense is valid,
where $[x]$ is an integer part of a real number $x;$
another notations are the same as in Theorem~{\rm 1}.}

\vspace{2mm}

It should be noted that an analogue of Theorem 2 was considered 
in \cite{Rybakov1000}. 
Note that we use another notations 
\cite{2018a} (Sect.~1.11), \cite{arxiv-1} (Sect.~15)
in comparison with \cite{Rybakov1000}.
Moreover, the proof of an analogue of Theorem 2
from \cite{Rybakov1000} is somewhat different from the proof given in 
\cite{2018a} (Sect.~1.11), \cite{arxiv-1} (Sect.~15).

\vspace{5mm}

\section{Calculation of the Mean-Square Approximation Error
in the Method of Generalized Multiple Fourier Seires}

\vspace{5mm}

Note that for the integrals $J[\psi^{(k)}]_{T,t}$ defined by 
(\ref{ito})
the mean-square approximation error can be exactly
calculated and efficiently estimated.

Let $J[\psi^{(k)}]_{T,t}^{q}$ be the
expression on the right-hand side of (\ref{leto6000}) before passing to the limit
$\hbox{\vtop{\offinterlineskip\halign{
\hfil#\hfil\cr
{\rm l.i.m.}\cr
$\stackrel{}{{}_{p_1,\ldots,p_k\to \infty}}$\cr
}} }$ for the case
$p_1=\ldots=p_k=q,$ i.e.

\vspace{1mm}

$$
J[\psi^{(k)}]_{T,t}^{q}=
\sum\limits_{j_1,\ldots,j_k=0}^{q}
C_{j_k\ldots j_1}\Biggl(
\prod_{l=1}^k\zeta_{j_l}^{(i_l)}+\sum\limits_{r=1}^{[k/2]}
(-1)^r \times
\Biggr.
$$

\vspace{2mm}
\begin{equation}
\label{r1}
\times
\sum_{\stackrel{(\{\{g_1, g_2\}, \ldots, 
\{g_{2r-1}, g_{2r}\}\}, \{q_1, \ldots, q_{k-2r}\})}
{{}_{\{g_1, g_2, \ldots, 
g_{2r-1}, g_{2r}, q_1, \ldots, q_{k-2r}\}=\{1, 2, \ldots, k\}}}}
\prod\limits_{s=1}^r
{\bf 1}_{\{i_{g_{{}_{2s-1}}}=~i_{g_{{}_{2s}}}\ne 0\}}
\Biggl.{\bf 1}_{\{j_{g_{{}_{2s-1}}}=~j_{g_{{}_{2s}}}\}}
\prod_{l=1}^{k-2r}\zeta_{j_{q_l}}^{(i_{q_l})}\Biggr)
\end{equation}

\vspace{5mm}

Let us denote

$$
{\sf M}\left\{\left(J[\psi^{(k)}]_{T,t}-
J[\psi^{(k)}]_{T,t}^{q}\right)^2\right\}\stackrel{{\rm def}}
{=}E_k^{q},
$$

\vspace{2mm}
$$
\int\limits_{[t,T]^k}
K^2(t_1,\ldots,t_k)dt_1\ldots dt_k
\stackrel{{\rm def}}{=}I_k.
$$

\vspace{4mm}

In \cite{2017-1}-\cite{2018aa1}, \cite{arxiv-1}, 
\cite{arxiv-2} it was shown that 

\vspace{-1mm}
\begin{equation}
\label{qq4}
E_k^{q}\le k!\left(I_k-\sum_{j_1,\ldots,j_k=0}^{q}C^2_{j_k\ldots j_1}\right)
\end{equation}

\vspace{2mm}
\noindent
for the following two cases:

\vspace{2mm}

1.\ $i_1,\ldots,i_k=1,\ldots,m$ and $T-t\in (0, +\infty)$,

2.\ $i_1,\ldots,i_k=0, 1,\ldots,m$ and  $T-t\in (0, 1)$.

\vspace{3mm}

The value $E_k^{q}$
can be calculated exactly.

\vspace{2mm}

{\bf Theorem 3} \cite{2018a} (Sect.~1.12), \cite{arxiv-2} (Sect.~6).
{\it Suppose that $\{\phi_j(x)\}_{j=0}^{\infty}$ 
is an arbitrary complete orthonormal system  
of functions in the space $L_2([t,T])$ and
$\psi_1(\tau),\ldots,\psi_k(\tau)\in L_2([t, T]).$  
Then

\begin{equation}
\label{tttr11}
E_k^q=I_k- \sum_{j_1,\ldots, j_k=0}^{q}
C_{j_k\ldots j_1}
{\sf M}\left\{J[\psi^{(k)}]_{T,t}
\sum\limits_{(j_1,\ldots,j_k)}
\int\limits_t^T \phi_{j_k}(t_k)
\ldots
\int\limits_t^{t_{2}}\phi_{j_{1}}(t_{1})
d{\bf f}_{t_1}^{(i_1)}\ldots
d{\bf f}_{t_k}^{(i_k)}\right\},
\end{equation}

\vspace{5mm}
\noindent
where
$i_1,\ldots,i_k = 1,\ldots,m;$\ 
the expression 

\vspace{-1mm}
$$
\sum\limits_{(j_1,\ldots,j_k)}
$$ 

\vspace{3mm}
\noindent
means the sum with respect to all
possible permutations
$(j_1,\ldots,j_k)$. At the same time if 
$j_r$ swapped with $j_q$ in the permutation $(j_1,\ldots,j_k)$,
then $i_r$ swapped with $i_q$ in the permutation 
$(i_1,\ldots,i_k);$
another notations are the same as in Theorems {\rm 1, 2.}
}

\vspace{2mm}

Note that 

\vspace{-2mm}
$$
{\sf M}\left\{J[\psi^{(k)}]_{T,t}
\int\limits_t^T \phi_{j_k}(t_k)
\ldots
\int\limits_t^{t_{2}}\phi_{j_{1}}(t_{1})
d{\bf f}_{t_1}^{(i_1)}\ldots
d{\bf f}_{t_k}^{(i_k)}\right\}=C_{j_k\ldots j_1}.
$$

\vspace{4mm}

Therefore, for the special case of pairwise 
different numbers $i_1,\ldots,i_k$ 
as well as for the case $i_1=\ldots=i_k$
from Theorem 3 it follows that
\cite{2018}-\cite{2018aa1}, \cite{17a}, 
\cite{arxiv-2}

\vspace{-1mm}
\begin{equation}
\label{qq1}
E_k^q= I_k- \sum_{j_1,\ldots,j_k=0}^{q}
C_{j_k\ldots j_1}^2,
\end{equation}

\vspace{2mm}
$$
E_k^q= I_k - \sum_{j_1,\ldots,j_k=0}^{q}
C_{j_k\ldots j_1}\Biggl(\sum\limits_{(j_1,\ldots,j_k)}
C_{j_k\ldots j_1}\Biggr),
$$

\vspace{3mm}
\noindent
where 

\vspace{-2mm}
$$
\sum\limits_{(j_1,\ldots,j_k)}
$$ 

\vspace{2mm}
\noindent
is a sum with respect to all 
possible permutations
$(j_1,\ldots,j_k)$.

Consider some examples \cite{2018}-\cite{2018aa1}, \cite{17a}, 
\cite{arxiv-2} of application of Theorem 3
$(i_1,i_2,i_3=1,\ldots,m)$

\vspace{1mm}
\begin{equation}
\label{qq2}
E_2^q     
=I_2
-\sum_{j_1,j_2=0}^q
C_{j_2j_1}^2-
\sum_{j_1,j_2=0}^q
C_{j_2j_1}C_{j_1j_2}\ \ \ (i_1=i_2),
\end{equation}

\vspace{1mm}
\begin{equation}
\label{qq3}
E_3^q=I_3
-\sum_{j_3,j_2,j_1=0}^q C_{j_3j_2j_1}^2-
\sum_{j_3,j_2,j_1=0}^q C_{j_3j_1j_2}C_{j_3j_2j_1}\ \ \ (i_1=i_2\ne i_3),
\end{equation}

\vspace{1mm}
\begin{equation}
\label{882}
E_3^q=I_3-
\sum_{j_3,j_2,j_1=0}^q C_{j_3j_2j_1}^2-
\sum_{j_3,j_2,j_1=0}^q C_{j_2j_3j_1}C_{j_3j_2j_1}\ \ \ (i_1\ne i_2=i_3),
\end{equation}

\vspace{1mm}
\begin{equation}
\label{883}
E_3^q=I_3
-\sum_{j_3,j_2,j_1=0}^q C_{j_3j_2j_1}^2-
\sum_{j_3,j_2,j_1=0}^q C_{j_3j_2j_1}C_{j_1j_2j_3}\ \ \ (i_1=i_3\ne i_2).
\end{equation}

\vspace{6mm}

The values $E_4^q$ and $E_5^q$ were calculated exaclty for all possible 
combinations of $i_1,\ldots,i_5=1,\ldots,m$ in 
\cite{2018}-\cite{2018aa1},  
\cite{arxiv-2}.

\vspace{5mm}

\section{Expansions 
of Iterated Stratonovich
Stochastic Integrals Based on Multiple Fourier--Legendre
Series and Multiple Trigonometric Fourier Seires}

\vspace{5mm}

In contrast to the iterated Ito stochastic integrals, 
the iterated Stratonovich stochastic integrals
have simpler expansions than (\ref{tyyy}) and (\ref{leto6000}), but the calculation (or estimation) 
of the mean-square approximation
errors for the latter is a much more difficult problem than 
for the former. Study this issue in
detail.

Introduce the following 
iterated Stratonovich stochastic integrals 

\vspace{-1mm}
\begin{equation}
\label{str}
J^{*}[\psi^{(k)}]_{T,t}=
{\int\limits_t^{*}}^T
\psi_k(t_k) \ldots 
{\int\limits_t^{*}}^{t_2}
\psi_1(t_1) d{\bf w}_{t_1}^{(i_1)}\ldots
d{\bf w}_{t_k}^{(i_k)},
\end{equation}

\vspace{3mm}
\noindent
where the notations are the same as in the formula (\ref{ito})
(in this paper, 
we use the definition of the Stratonovich stochastic integral from \cite{KlPl2}).

Consider a slightly modified and extended theoretical result 
that adapts Theorems 1, 2 for
the iterated Stratonovich stochastic
integrals (\ref{str}) of multiplicities 2 to 4 (some old results).

\vspace{2mm}

{\bf Theorem 4}\ \cite{2011-2}-\cite{2018aa1}, \cite{2010-2}-\cite{2013},
\cite{30a}, \cite{271a}, \cite{arxiv-3}-\cite{arxiv-8}, \cite{arxiv-6}.\ {\it 
Assume that $\{\phi_j(x)\}_{j=0}^{\infty}$ 
is a complete or\-tho\-nor\-mal system of Legendre
polynomials or trigonometric functions in $L_2([t, T]).$ 
In addition, assume that $\psi_2(s)$ is a continuously
differentiable function on the interval $[t, T]$  
and $\psi_1(s),$ $\psi_3(s)$ are twice continuous\-ly differentiable
functions on the interval $[t, T]$. Then

\vspace{-1mm}
\begin{equation}
\label{zz}
J^{*}[\psi^{(k)}]_{T,t}=
\hbox{\vtop{\offinterlineskip\halign{
\hfil#\hfil\cr
{\rm l.i.m.}\cr
$\stackrel{}{{}_{q\to \infty}}$\cr
}} }
\sum\limits_{j_1, \ldots, j_k=0}^{q}
C_{j_k \ldots j_1}\zeta_{j_1}^{(i_1)} \ldots \zeta_{j_k}^{(i_k)},
\end{equation}

\vspace{3mm}
\noindent
where $k=2, 3, 4.$ At that $\psi_1(s),\ldots, \psi_k(s)\equiv 1$ and
$i_1,\ldots,i_k=0, 1,\ldots,m$
in {\rm (\ref{zz})} 
for $k=4,$ while $i_1,\ldots,i_k=1,\ldots,m$ in {\rm (\ref{zz})} 
for $k=2, 3;$
the other notations are the same as in Theorems {\rm 1, 2}.
}

\vspace{2mm}

Recently, a new approach to the expansion and mean-square 
approximation of iterated Stratonovich stochastic integrals has been obtained
\cite{2018a} (Sect.~2.10--2.16), \cite{arxiv-4} (Sect.~7--13), \cite{arxiv-5} (Sect.~13--19), 
\cite{arxiv-6} (Sect.~5--11), \cite{new-art-1-xxy}
(Sect.~4--9).
Let us formulate four theorems that were obtained using this approach.

\vspace{2mm}  

{\bf Theorem 5}\ \cite{2018a}, \cite{arxiv-4}, \cite{arxiv-5}, 
\cite{arxiv-6}, \cite{new-art-1-xxy}.\
{\it Suppose 
that $\{\phi_j(x)\}_{j=0}^{\infty}$ is a complete orthonormal system of 
Legendre polynomials or trigonometric functions in the space $L_2([t, T]).$
Furthermore, let $\psi_1(\tau), \psi_2(\tau),$ $\psi_3(\tau)$ are continuously dif\-ferentiable 
nonrandom functions on $[t, T].$ 
Then, for the 
iterated Stra\-to\-no\-vich stochastic integral of third multiplicity

$$
J^{*}[\psi^{(3)}]_{T,t}={\int\limits_t^{*}}^T\psi_3(t_3)
{\int\limits_t^{*}}^{t_3}\psi_2(t_2)
{\int\limits_t^{*}}^{t_2}\psi_1(t_1)
d{\bf w}_{t_1}^{(i_1)}
d{\bf w}_{t_2}^{(i_2)}d{\bf w}_{t_3}^{(i_3)}\ \ \ (i_1,i_2,i_3=0,1,\ldots,m)
$$

\vspace{4mm}
\noindent
the following 
relations

\vspace{-1mm}
\begin{equation}
\label{fin1}
J^{*}[\psi^{(3)}]_{T,t}
=\hbox{\vtop{\offinterlineskip\halign{
\hfil#\hfil\cr
{\rm l.i.m.}\cr
$\stackrel{}{{}_{p\to \infty}}$\cr
}} }
\sum\limits_{j_1, j_2, j_3=0}^{p}
C_{j_3 j_2 j_1}\zeta_{j_1}^{(i_1)}\zeta_{j_2}^{(i_2)}\zeta_{j_3}^{(i_3)},
\end{equation}

\vspace{3mm}
\begin{equation}
\label{fin2}
{\sf M}\left\{\left(
J^{*}[\psi^{(3)}]_{T,t}-
\sum\limits_{j_1, j_2, j_3=0}^{p}
C_{j_3 j_2 j_1}\zeta_{j_1}^{(i_1)}\zeta_{j_2}^{(i_2)}\zeta_{j_3}^{(i_3)}\right)^2\right\}
\le \frac{C}{p}
\end{equation}

\vspace{5mm}
\noindent
are fulfilled, where $i_1, i_2, i_3=0,1,\ldots,m$ in {\rm (\ref{fin1})} and 
$i_1, i_2, i_3=1,\ldots,m$ in {\rm (\ref{fin2}),}
constant $C$ is independent of $p,$

$$
C_{j_3 j_2 j_1}=\int\limits_t^T\psi_3(t_3)\phi_{j_3}(t_3)
\int\limits_t^{t_3}\psi_2(t_2)\phi_{j_2}(t_2)
\int\limits_t^{t_2}\psi_1(t_1)\phi_{j_1}(t_1)dt_1dt_2dt_3
$$

\vspace{4mm}
\noindent
and
$$
\zeta_{j}^{(i)}=
\int\limits_t^T \phi_{j}(\tau) d{\bf f}_{\tau}^{(i)}
$$ 

\vspace{2mm}
\noindent
are independent standard Gaussian random variables for various 
$i$ or $j$ {\rm (}in the case when $i\ne 0${\rm );} 
another notations are the same as in Theorems~{\rm 1, 2}.}

\vspace{2mm}

{\bf Theorem 6}\ \cite{2018a}, \cite{arxiv-4}, \cite{arxiv-5}, 
\cite{arxiv-6}, \cite{new-art-1-xxy}.\ {\it Let
$\{\phi_j(x)\}_{j=0}^{\infty}$ be a complete orthonormal system of 
Legendre polynomials or trigonometric functions in the space $L_2([t, T]).$
Furthermore, let $\psi_1(\tau), \ldots, \psi_4(\tau)$ be continuously dif\-ferentiable 
nonrandom functions on $[t, T].$ 
Then, for the 
iterated Stra\-to\-no\-vich stochastic integral of fourth multiplicity

\begin{equation}
\label{fin0}
J^{*}[\psi^{(4)}]_{T,t}={\int\limits_t^{*}}^T\psi_4(t_4)
{\int\limits_t^{*}}^{t_4}\psi_3(t_3)
{\int\limits_t^{*}}^{t_3}\psi_2(t_2)
{\int\limits_t^{*}}^{t_2}\psi_1(t_1)
d{\bf w}_{t_1}^{(i_1)}
d{\bf w}_{t_2}^{(i_2)}d{\bf w}_{t_3}^{(i_3)}d{\bf w}_{t_4}^{(i_4)}
\end{equation}

\vspace{4mm}
\noindent
the following 
relations

\begin{equation}
\label{fin3}
J^{*}[\psi^{(4)}]_{T,t}
=\hbox{\vtop{\offinterlineskip\halign{
\hfil#\hfil\cr
{\rm l.i.m.}\cr
$\stackrel{}{{}_{p\to \infty}}$\cr
}} }
\sum\limits_{j_1, j_2, j_3,j_4=0}^{p}
C_{j_4j_3 j_2 j_1}\zeta_{j_1}^{(i_1)}\zeta_{j_2}^{(i_2)}\zeta_{j_3}^{(i_3)}\zeta_{j_4}^{(i_4)},
\end{equation}

\vspace{3mm}

\begin{equation}
\label{fin4}
{\sf M}\left\{\left(
J^{*}[\psi^{(4)}]_{T,t}-
\sum\limits_{j_1, j_2, j_3, j_4=0}^{p}
C_{j_4 j_3 j_2 j_1}\zeta_{j_1}^{(i_1)}\zeta_{j_2}^{(i_2)}\zeta_{j_3}^{(i_3)}
\zeta_{j_4}^{(i_4)}
\right)^2\right\}
\le \frac{C}{p^{1-\varepsilon}}
\end{equation}

\vspace{5mm}
\noindent
are fulfilled, where $i_1, \ldots , i_4=0,1,\ldots,m$ in {\rm (\ref{fin0}),} {\rm (\ref{fin3})} 
and $i_1, \ldots, i_4=1,\ldots,m$ in {\rm (\ref{fin4}),}
constant $C$ does not depend on $p,$
$\varepsilon$ is an arbitrary
small positive real number 
for the case of complete orthonormal system of 
Legendre polynomials in the space $L_2([t, T])$
and $\varepsilon=0$ for the case of
complete orthonormal system of 
trigonometric functions in the space $L_2([t, T]),$

$$
C_{j_4 j_3 j_2 j_1}=
$$

$$
=
\int\limits_t^T\psi_4(t_4)\phi_{j_4}(t_4)
\int\limits_t^{t_4}\psi_3(t_3)\phi_{j_3}(t_3)
\int\limits_t^{t_3}\psi_2(t_2)\phi_{j_2}(t_2)
\int\limits_t^{t_2}\psi_1(t_1)\phi_{j_1}(t_1)dt_1dt_2dt_3dt_4;
$$

\vspace{4mm}
\noindent
another notations are the same as in Theorem~{\rm 5}.}

\vspace{2mm}

{\bf Theorem 7}\ \cite{2018a}, \cite{arxiv-4}, \cite{arxiv-5}, 
\cite{arxiv-6}, \cite{new-art-1-xxy}.\
{\it Assume 
that $\{\phi_j(x)\}_{j=0}^{\infty}$ is a complete orthonormal system of 
Legendre polynomials or trigonometric functions in the space $L_2([t, T])$
and $\psi_1(\tau), \ldots, \psi_5(\tau)$ are continuously dif\-ferentiable 
nonrandom functions on $[t, T].$ 
Then, for the 
iterated Stra\-to\-no\-vich stochastic integral of fifth multiplicity

\begin{equation}
\label{fin7}
J^{*}[\psi^{(5)}]_{T,t}={\int\limits_t^{*}}^T\psi_5(t_5)
\ldots
{\int\limits_t^{*}}^{t_2}\psi_1(t_1)
d{\bf w}_{t_1}^{(i_1)}
\ldots d{\bf w}_{t_5}^{(i_5)}
\end{equation}

\vspace{4mm}
\noindent
the following 
relations

\begin{equation}
\label{fin8}
J^{*}[\psi^{(5)}]_{T,t}
=\hbox{\vtop{\offinterlineskip\halign{
\hfil#\hfil\cr
{\rm l.i.m.}\cr
$\stackrel{}{{}_{p\to \infty}}$\cr
}} }
\sum\limits_{j_1,\ldots,j_5=0}^{p}
C_{j_5 \ldots j_1}\zeta_{j_1}^{(i_1)}\ldots \zeta_{j_5}^{(i_5)},
\end{equation}

\vspace{3mm}

\begin{equation}
\label{fin9}
{\sf M}\left\{\left(
J^{*}[\psi^{(5)}]_{T,t}-
\sum\limits_{j_1, \ldots, j_5=0}^{p}
C_{j_5 \ldots j_1}\zeta_{j_1}^{(i_1)}\ldots
\zeta_{j_5}^{(i_5)}
\right)^2\right\}
\le \frac{C}{p^{1-\varepsilon}}
\end{equation}

\vspace{5mm}
\noindent
are fulfilled, where $i_1, \ldots , i_5=0,1,\ldots,m$ in {\rm (\ref{fin7}),} {\rm (\ref{fin8})} 
and $i_1, \ldots, i_5=1,\ldots,m$ in {\rm (\ref{fin9}),}
constant $C$ is independent of $p,$
$\varepsilon$ is an arbitrary
small positive real number 
for the case of complete orthonormal system of 
Legendre polynomials in the space $L_2([t, T])$
and $\varepsilon=0$ for the case of
complete orthonormal system of 
trigonometric functions in the space $L_2([t, T]),$

$$
C_{j_5 \ldots j_1}=
\int\limits_t^T\psi_5(t_5)\phi_{j_5}(t_5)\ldots
\int\limits_t^{t_2}\psi_1(t_1)\phi_{j_1}(t_1)dt_1\ldots dt_5;
$$

\vspace{3mm}
\noindent
another notations are the same as in Theorems~{\rm 5, 6}.}

\vspace{2mm}

{\bf Theorem 8}\ \cite{2018a}, \cite{arxiv-4}, \cite{arxiv-5}, 
\cite{arxiv-6}.\
{\it Suppose that 
$\{\phi_j(x)\}_{j=0}^{\infty}$ is a complete orthonormal system of 
Legendre polynomials or trigonometric functions in the space $L_2([t, T]).$
Then, for the 
iterated Stratonovich stochastic integral of sixth multiplicity

\begin{equation}
\label{after10001qu1}
J_{T,t}^{*(i_1\ldots i_6)}={\int\limits_t^{*}}^T
\ldots
{\int\limits_t^{*}}^{t_2}
d{\bf w}_{t_1}^{(i_1)}
\ldots d{\bf w}_{t_6}^{(i_6)}
\end{equation}

\vspace{3mm}
\noindent
the following 
expansion 

\vspace{-1mm}
$$
J_{T,t}^{*(i_1\ldots i_6)}
=\hbox{\vtop{\offinterlineskip\halign{
\hfil#\hfil\cr
{\rm l.i.m.}\cr
$\stackrel{}{{}_{p\to \infty}}$\cr
}} }
\sum\limits_{j_1, \ldots, j_6=0}^{p}
C_{j_6 \ldots j_1}\zeta_{j_1}^{(i_1)}\ldots
\zeta_{j_6}^{(i_6)}
$$

\vspace{5mm}
\noindent
that converges in the mean-square sense is valid, where
$i_1, \ldots, i_6=0, 1,\ldots,m,$

$$
C_{j_6 \ldots j_1}=
\int\limits_t^T\phi_{j_6}(t_6)\ldots
\int\limits_t^{t_2}\phi_{j_1}(t_1)dt_1\ldots dt_6;
$$

\vspace{3mm}
\noindent
another notations are the same as in Theorems~{\rm 5--7}.}

\vspace{2mm}

Clearly, the expansion (\ref{zz}) is simpler than the expansions (\ref{tyyy}), (\ref{leto6000}).
However, the calculation of the mean-square
approximation error for the expansion (\ref{zz}) turns out to be 
much more difficult than for the expansions
(\ref{tyyy}), (\ref{leto6000}). We will demonstrate this fact below.

The cases $k=1, 2$ are actually not interesting: for $k=1$, 
the Ito and Stratonovich stochastic
integrals of a smooth nonrandom function equal each other with 
probability 1 (w.~p.~1); for $k=2$, the
Ito stochastic integrals appearing in the numerical scheme 
(\ref{4.45}) differ w.~p.~1 from the corresponding
Stratonovich stochastic integrals by constant values by virtue of the 
standard relations between the
Ito and Stratonovich stochastic integrals \cite{KlPl2}.
Consider the triple Stratonovich stochastic integral defined by

\vspace{-1mm}
$$
I_{(000)T,t}^{*(i_1i_2i_3)}=
\int\limits_t^{*T}\int\limits_t^{*t_{3}}
\int\limits_t^{*t_{2}}
d{\bf f}_{t_1}^{(i_1)}
d{\bf f}_{t_2}^{(i_2)}
d{\bf f}_{t_3}^{(i_3)}\ \ \ (i_1, i_2, i_3=1,\ldots,m).
$$

\vspace{3mm}

In view of the standard relations between 
the Ito and Stratonovich stochastic integrals \cite{KlPl2} and
also Theorems 1, 2, and 4 ($k=3$), we obtain

$$
{\sf M}\left\{\left(I_{(000)T,t}^{*(i_1i_2i_3)}-
I_{(000)T,t}^{*(i_1i_2i_3)q}\right)^2\right\}=
$$

\vspace{2mm}
$$
={\sf M}\left\{\left(I_{(000)T,t}^{(i_1i_2i_3)}+
{\bf 1}_{\{i_1=i_2\}}
\frac{1}{2}\int\limits_t^T\int\limits_t^{\tau}dsd{\bf f}_{\tau}^{(i_3)}
+{\bf 1}_{\{i_2=i_3\}}
\frac{1}{2}\int\limits_t^T\int\limits_t^{\tau}d{\bf f}_{s}^{(i_1)}d\tau-
I_{(000)T,t}^{*(i_1i_2i_3)q}\right)^2\right\}=
$$

\vspace{2mm}
$$
={\sf M}\Biggl\{\Biggl(I_{(000)T,t}^{(i_1i_2i_3)}-I_{(000)T,t}^{(i_1i_2i_3)q}+
I_{(000)T,t}^{(i_1i_2i_3)q}+
\Biggr.\Biggr.
$$

\vspace{2mm}
\begin{equation}
\label{tango3}
+
\left.\Biggl.
{\bf 1}_{\{i_1=i_2\}}
\frac{1}{2}\int\limits_t^T\int\limits_t^{\tau}dsd{\bf f}_{\tau}^{(i_3)}
+{\bf 1}_{\{i_2=i_3\}}
\frac{1}{2}\int\limits_t^T\int\limits_t^{\tau}d{\bf f}_{s}^{(i_1)}d\tau-
I_{(000)T,t}^{*(i_1i_2i_3)q}\Biggr)^2\right\},
\end{equation}

\vspace{6mm}

$$
I_{(000)T,t}^{(i_1i_2i_3)q}=
\sum_{j_1,j_2,j_3=0}^{q}
C_{j_3j_2j_1}\Biggl(
\zeta_{j_1}^{(i_1)}\zeta_{j_2}^{(i_2)}\zeta_{j_3}^{(i_3)}
-{\bf 1}_{\{i_1=i_2\}}
{\bf 1}_{\{j_1=j_2\}}
\zeta_{j_3}^{(i_3)}-
\Biggr.
$$

\vspace{2mm}
\begin{equation}
\label{tango1}
\Biggl.
-{\bf 1}_{\{i_2=i_3\}}
{\bf 1}_{\{j_2=j_3\}}
\zeta_{j_1}^{(i_1)}-
{\bf 1}_{\{i_1=i_3\}}
{\bf 1}_{\{j_1=j_3\}}
\zeta_{j_2}^{(i_2)}\Biggr),
\end{equation}

\vspace{5mm}

\begin{equation}
\label{tango2}
I_{(000)T,t}^{*(i_1i_2i_3)q}=
\sum_{j_1,j_2,j_3=0}^{q}
C_{j_3j_2j_1}
\zeta_{j_1}^{(i_1)}\zeta_{j_2}^{(i_2)}\zeta_{j_3}^{(i_3)},
\end{equation}

\vspace{6mm}
\noindent
where $I_{(000)T,t}^{(i_1i_2i_3)q}$ 
is the approximation defined by (\ref{r1})
(also see (\ref{a3})) 
for $k=3$
and $I_{(000)T,t}^{*(i_1i_2i_3)q}$ 
is the approximation
based on Theorem 4 for
$k=3.$

Substituting (\ref{tango1}) and (\ref{tango2}) into (\ref{tango3}) yields

\vspace{2mm}
$$
{\sf M}\left\{\left(I_{(000)T,t}^{*(i_1i_2i_3)}-
I_{(000)T,t}^{*(i_1i_2i_3)q}\right)^2\right\}=
$$

\vspace{2mm}
$$
={\sf M}\left\{\left(
I_{(000)T,t}^{(i_1i_2i_3)}-I_{(000)T,t}^{(i_1i_2i_3)q}+
{\bf 1}_{\{i_1=i_2\}}
\left(
\frac{1}{2}\int\limits_t^T\int\limits_t^{\tau}dsd{\bf f}_{\tau}^{(i_3)}-
\sum_{j_1,j_3=0}^{q}
C_{j_3j_1j_1}
\zeta_{j_3}^{(i_3)}\right)+\right.\right.
$$

\vspace{2mm}
$$
+{\bf 1}_{\{i_2=i_3\}}\left(
\frac{1}{2}\int\limits_t^T\int\limits_t^{\tau}d{\bf f}_{s}^{(i_1)}d\tau-
\sum_{j_1,j_3=0}^{q}
C_{j_3j_3j_1}
\zeta_{j_1}^{(i_1)}\right)
\left.\left.-{\bf 1}_{\{i_1=i_3\}}
\sum_{j_1,j_2=0}^{q}
C_{j_1j_2j_1}
\zeta_{j_2}^{(i_2)}\right)^2\right\}\le
$$

\vspace{5mm}
$$
\le 4 \Biggl({\sf M}\left\{\left(I_{(000)T,t}^{(i_1i_2i_3)}-
I_{(000)T,t}^{(i_1i_2i_3)q}\right)^2\right\}+\Biggr.
{\bf 1}_{\{i_1=i_2\}}F^{(i_3)}_q+
$$

\vspace{2mm}
\begin{equation}
\label{tango4}
\Biggl.+
{\bf 1}_{\{i_2=i_3\}}G^{(i_1)}_q+
{\bf 1}_{\{i_1=i_3\}}H^{(i_2)}_q\Biggr),
\end{equation}

\vspace{2mm}
\noindent
where

$$
F^{(i_3)}_q={\sf M}\left\{\left(
\frac{1}{2}\int\limits_t^T\int\limits_t^{\tau}dsd{\bf f}_{\tau}^{(i_3)}-
\sum_{j_1,j_3=0}^{q}
C_{j_3j_1j_1}
\zeta_{j_3}^{(i_3)}\right)^2\right\},
$$

\vspace{2mm}
$$
G^{(i_1)}_q={\sf M}\left\{\left(
\frac{1}{2}\int\limits_t^T\int\limits_t^{\tau}d{\bf f}_{s}^{(i_1)}d\tau-
\sum_{j_1,j_3=0}^{q}
C_{j_3j_3j_1}
\zeta_{j_1}^{(i_1)}\right)^2\right\},
$$

\vspace{2mm}
$$
H^{(i_2)}_q={\sf M}\left\{\left(
\sum_{j_1,j_2=0}^{q}
C_{j_1j_2j_1}
\zeta_{j_2}^{(i_2)}\right)^2\right\}.
$$

\vspace{5mm}

For the case of Legendre polynomials or
trigonometric functions we have the equalities
(for details, see the proof of Theorem 4 for $k=3$ in 
\cite{2011-2}-\cite{2018aa1}, \cite{2010-2}-\cite{2013},
\cite{arxiv-5})

\vspace{2mm}
$$
\lim\limits_{q\to\infty}F^{(i_3)}_q=0,\ \ \
\lim\limits_{q\to\infty} G^{(i_1)}_q=0,\ \ \ 
\lim\limits_{q\to\infty}H^{(i_2)}_q=0.
$$

\vspace{5mm}

However, in accordance with (\ref{tango4}) the value

\vspace{1mm}
$$
{\sf M}\left\{\left(I_{(000)T,t}^{*(i_1i_2i_3)}-
I_{(000)T,t}^{*(i_1i_2i_3)q}\right)^2\right\}
$$

\vspace{4mm}
\noindent
with a finite $q$ can be estimated by the sum of 

\vspace{1mm}
\begin{equation}
\label{r2}
4{\sf M}\left\{\left(I_{(000)T,t}^{(i_1i_2i_3)}-
I_{(000)T,t}^{(i_1i_2i_3)q}\right)^2\right\},
\end{equation}                                              

\vspace{4mm}
\noindent
and three additional terms of a rather complex structure. 
The value (\ref{r2}) can be calculated exactly
using Theorem 3 or estimated using (\ref{qq4}) for the case $k=3$.

As is easily observed, this peculiarity will also apply to 
the iterated Stratonovich stochastic
integrals of multiplicities 4 and 5, with the only 
difference that the number of additional terms like
$F^{(i_3)}_q$,
$G^{(i_1)}_q$, and $H^{(i_2)}_q$
will be considerably higher and their structure will be more complicated.
Therefore, the payment for a relatively simple approximation of 
the iterated Stratonovich
stochastic integrals (Theorems 4--8) in comparison with the iterated
Ito stochastic integrals (Theorems 1, 2) is a much more 
difficult calculation or estimation 
procedure of their mean-square
approximation errors (see Chapter 5 in \cite{2018a} for detail).
This well explains why the main emphasis of the paper is on the 
approximation of the iterated
Ito stochastic integrals figuring in the numerical scheme (\ref{4.45}). 
Their approximation involves Theorems
1, 2 for $k=1,\ldots,5$ and also a complete orthonormal system of 
Legendre polynomials in the
space $L_2([t,T])$. As was established in \cite{2018a}-\cite{2018aa1},
\cite{29a}, \cite{301a}, 
the Legendre polynomials 
have a series of advantages over the
trigonometric functions for the approximation of iterated stochastic 
integrals using Theorems 1, 2.

\vspace{5mm}

\section{Approximation of Iterated 
Ito Stochastic Integrals Based on Multiple
Fourier--Legendre Series}

\vspace{5mm}

Consider the approximations of the iterated Ito stochastic integrals 
that appear in the numerical
scheme (\ref{4.45}) using Theorems 1, 2 and the complete orthonormal system of 
Legendre polynomials in the
space $L_2([\tau_p,\tau_{p+1}])$
($\tau_p=p\Delta,$ $N\Delta= T,$ 
$p=0,1,\ldots,N$) \cite{2006}
(also see \cite{2011-2}-\cite{arxiv-24}, \cite{Mikh-1}-\cite{Mikh-2})

\vspace{2mm}
\begin{equation}
\label{yyy1aaa}
I_{(0)\tau_{p+1},\tau_p}^{(i_1)}=\sqrt{\Delta}\zeta_0^{(i_1)},
\end{equation}

\vspace{2mm}
\begin{equation}
\label{qqqq1}
I_{(00)\tau_{p+1},\tau_p}^{(i_1 i_2)q}=
\frac{\Delta}{2}\left(\zeta_0^{(i_1)}\zeta_0^{(i_2)}+\sum_{i=1}^{q}
\frac{1}{\sqrt{4i^2-1}}\left(
\zeta_{i-1}^{(i_1)}\zeta_{i}^{(i_2)}-
\zeta_i^{(i_1)}\zeta_{i-1}^{(i_2)}\right) - {\bf 1}_{\{i_1=i_2\}}\right),
\end{equation}

\vspace{5mm}

\begin{equation}
\label{yyy2}
I_{(1)\tau_{p+1},\tau_p}^{(i_1)}=
-\frac{{\Delta}^{3/2}}{2}\left(\zeta_0^{(i_1)}+
\frac{1}{\sqrt{3}}\zeta_1^{(i_1)}\right),
\end{equation}

\vspace{5mm}

$$
I_{(000)\tau_{p+1},\tau_p}^{(i_1i_2i_3)q}
=\sum_{j_1,j_2,j_3=0}^{q}
C_{j_3j_2j_1}\Biggl(
\zeta_{j_1}^{(i_1)}\zeta_{j_2}^{(i_2)}\zeta_{j_3}^{(i_3)}
-{\bf 1}_{\{i_1=i_2\}}
{\bf 1}_{\{j_1=j_2\}}
\zeta_{j_3}^{(i_3)}-
\Biggr.
$$
\begin{equation}
\label{yyy3}
\Biggl.
-{\bf 1}_{\{i_2=i_3\}}
{\bf 1}_{\{j_2=j_3\}}
\zeta_{j_1}^{(i_1)}-
{\bf 1}_{\{i_1=i_3\}}
{\bf 1}_{\{j_1=j_3\}}
\zeta_{j_2}^{(i_2)}\Biggr),
\end{equation}

\vspace{5mm}

$$
I_{(0000)\tau_{p+1},\tau_p}^{(i_1 i_2 i_3 i_4)q}
=\sum_{j_1,j_2,j_3,j_4=0}^{q}
C_{j_4 j_3 j_2 j_1}\Biggl(
\prod_{l=1}^4\zeta_{j_l}^{(i_l)}
-\Biggr.
$$

\vspace{-1mm}
$$
-
{\bf 1}_{\{i_1=i_2\}}
{\bf 1}_{\{j_1=j_2\}}
\zeta_{j_3}^{(i_3)}
\zeta_{j_4}^{(i_4)}
-
{\bf 1}_{\{i_1=i_3\}}
{\bf 1}_{\{j_1=j_3\}}
\zeta_{j_2}^{(i_2)}
\zeta_{j_4}^{(i_4)}-
$$

\vspace{-2.5mm}
$$
-
{\bf 1}_{\{i_1=i_4\}}
{\bf 1}_{\{j_1=j_4\}}
\zeta_{j_2}^{(i_2)}
\zeta_{j_3}^{(i_3)}
-
{\bf 1}_{\{i_2=i_3\}}
{\bf 1}_{\{j_2=j_3\}}
\zeta_{j_1}^{(i_1)}
\zeta_{j_4}^{(i_4)}-
$$

\vspace{-2.5mm}
$$
-
{\bf 1}_{\{i_2=i_4\}}
{\bf 1}_{\{j_2=j_4\}}
\zeta_{j_1}^{(i_1)}
\zeta_{j_3}^{(i_3)}
-
{\bf 1}_{\{i_3=i_4\}}
{\bf 1}_{\{j_3=j_4\}}
\zeta_{j_1}^{(i_1)}
\zeta_{j_2}^{(i_2)}+
$$

\vspace{-2.5mm}
$$
+
{\bf 1}_{\{i_1=i_2\}}
{\bf 1}_{\{j_1=j_2\}}
{\bf 1}_{\{i_3=i_4\}}
{\bf 1}_{\{j_3=j_4\}}+
$$

\vspace{-2.5mm}
$$
+
{\bf 1}_{\{i_1=i_3\}}
{\bf 1}_{\{j_1=j_3\}}
{\bf 1}_{\{i_2=i_4\}}
{\bf 1}_{\{j_2=j_4\}}+
$$
\begin{equation}
\label{yyy4}
+\Biggl.
{\bf 1}_{\{i_1=i_4\}}
{\bf 1}_{\{j_1=j_4\}}
{\bf 1}_{\{i_2=i_3\}}
{\bf 1}_{\{j_2=j_3\}}\Biggr),
\end{equation}

\vspace{7mm}

$$
I_{(01)\tau_{p+1},\tau_p}^{(i_1 i_2)q}=
-\frac{\Delta}{2}
I_{(00)\tau_{p+1},\tau_p}^{(i_1 i_2)q}
-\frac{{\Delta}^2}{4}\Biggl(
\frac{1}{\sqrt{3}}\zeta_0^{(i_1)}\zeta_1^{(i_2)}+\Biggr.
$$

\vspace{1mm}
\begin{equation}
\label{yyy5}
+\Biggl.\sum_{i=0}^{q}\Biggl(
\frac{(i+2)\zeta_i^{(i_1)}\zeta_{i+2}^{(i_2)}
-(i+1)\zeta_{i+2}^{(i_1)}\zeta_{i}^{(i_2)}}
{\sqrt{(2i+1)(2i+5)}(2i+3)}-
\frac{\zeta_i^{(i_1)}\zeta_{i}^{(i_2)}}{(2i-1)(2i+3)}\Biggr)\Biggr),
\end{equation}

\vspace{8mm}

$$
I_{(10)\tau_{p+1},\tau_p}^{(i_1 i_2)q}=
-\frac{\Delta}{2}I_{(00)\tau_{p+1},\tau_p}^{(i_1 i_2)q}
-\frac{\Delta^2}{4}\Biggl(
\frac{1}{\sqrt{3}}\zeta_0^{(i_2)}\zeta_1^{(i_1)}+\Biggr.
$$

\vspace{1mm}
\begin{equation}
\label{yyy6}
+\Biggl.\sum_{i=0}^{q}\Biggl(
\frac{(i+1)\zeta_{i+2}^{(i_2)}\zeta_{i}^{(i_1)}
-(i+2)\zeta_{i}^{(i_2)}\zeta_{i+2}^{(i_1)}}
{\sqrt{(2i+1)(2i+5)}(2i+3)}+
\frac{\zeta_i^{(i_1)}\zeta_{i}^{(i_2)}}{(2i-1)(2i+3)}\Biggr)\Biggr)
\end{equation}

\vspace{8mm}
\noindent
or
$$
I_{(01)\tau_{p+1},\tau_p}^{(i_1 i_2)q}=
\sum_{j_1,j_2=0}^{q}
C_{j_2j_1}^{01}\Biggl(\zeta_{j_1}^{(i_1)}\zeta_{j_2}^{(i_2)}
-{\bf 1}_{\{i_1=i_2\}}
{\bf 1}_{\{j_1=j_2\}}\Biggr),
$$

\vspace{3mm}
$$
I_{(10)\tau_{p+1},\tau_p}^{(i_1 i_2)q}=
\sum_{j_1,j_2=0}^{q}
C_{j_2j_1}^{10}\Biggl(\zeta_{j_1}^{(i_1)}\zeta_{j_2}^{(i_2)}
-{\bf 1}_{\{i_1=i_2\}}
{\bf 1}_{\{j_1=j_2\}}\Biggr);
$$

\vspace{7mm}

\begin{equation}
\label{zzzz1}
{I}_{(2)\tau_{p+1},\tau_p}^{(i_1)}=
\frac{\Delta^{5/2}}{3}\left(
\zeta_0^{(i_1)}+\frac{\sqrt{3}}{2}\zeta_1^{(i_1)}+
\frac{1}{2\sqrt{5}}\zeta_2^{(i_1)}\right),
\end{equation}

\vspace{6mm}

$$
I_{(001)\tau_{p+1},\tau_p}^{(i_1i_2i_3)q}
=\sum_{j_1,j_2,j_3=0}^{q}
C_{j_3j_2j_1}^{001}\Biggl(
\zeta_{j_1}^{(i_1)}\zeta_{j_2}^{(i_2)}\zeta_{j_3}^{(i_3)}
-{\bf 1}_{\{i_1=i_2\}}
{\bf 1}_{\{j_1=j_2\}}
\zeta_{j_3}^{(i_3)}-
\Biggr.
$$
\begin{equation}
\label{yyy7}
\Biggl.
-{\bf 1}_{\{i_2=i_3\}}
{\bf 1}_{\{j_2=j_3\}}
\zeta_{j_1}^{(i_1)}-
{\bf 1}_{\{i_1=i_3\}}
{\bf 1}_{\{j_1=j_3\}}
\zeta_{j_2}^{(i_2)}\Biggr),
\end{equation}

\vspace{6mm}

$$
I_{(010)\tau_{p+1},\tau_p}^{(i_1i_2i_3)q}
=\sum_{j_1,j_2,j_3=0}^{q}
C_{j_3j_2j_1}^{010}\Biggl(
\zeta_{j_1}^{(i_1)}\zeta_{j_2}^{(i_2)}\zeta_{j_3}^{(i_3)}
-{\bf 1}_{\{i_1=i_2\}}
{\bf 1}_{\{j_1=j_2\}}
\zeta_{j_3}^{(i_3)}-
\Biggr.
$$
\begin{equation}
\label{yyy8}
\Biggl.
-{\bf 1}_{\{i_2=i_3\}}
{\bf 1}_{\{j_2=j_3\}}
\zeta_{j_1}^{(i_1)}-
{\bf 1}_{\{i_1=i_3\}}
{\bf 1}_{\{j_1=j_3\}}
\zeta_{j_2}^{(i_2)}\Biggr),
\end{equation}

\vspace{6mm}

$$
I_{(100)\tau_{p+1},\tau_p}^{(i_1i_2i_3)q}
=\sum_{j_1,j_2,j_3=0}^{q}
C_{j_3j_2j_1}^{100}\Biggl(
\zeta_{j_1}^{(i_1)}\zeta_{j_2}^{(i_2)}\zeta_{j_3}^{(i_3)}
-{\bf 1}_{\{i_1=i_2\}}
{\bf 1}_{\{j_1=j_2\}}
\zeta_{j_3}^{(i_3)}-
\Biggr.
$$
\begin{equation}
\label{yyy9}
\Biggl.
-{\bf 1}_{\{i_2=i_3\}}
{\bf 1}_{\{j_2=j_3\}}
\zeta_{j_1}^{(i_1)}-
{\bf 1}_{\{i_1=i_3\}}
{\bf 1}_{\{j_1=j_3\}}
\zeta_{j_2}^{(i_2)}\Biggr),
\end{equation}

\vspace{6mm}

$$
I_{(00000)\tau_{p+1},\tau_p}^{(i_1 i_2 i_3 i_4 i_5)q}
=\sum_{j_1,j_2,j_3,j_4,j_5=0}^q
C_{j_5 j_4 j_3 j_2 j_1}\Biggl(
\prod_{l=1}^5\zeta_{j_l}^{(i_l)}
-\Biggr.
$$
$$
-
{\bf 1}_{\{j_1=j_2\}}
{\bf 1}_{\{i_1=i_2\}}
\zeta_{j_3}^{(i_3)}
\zeta_{j_4}^{(i_4)}
\zeta_{j_5}^{(i_5)}-
{\bf 1}_{\{j_1=j_3\}}
{\bf 1}_{\{i_1=i_3\}}
\zeta_{j_2}^{(i_2)}
\zeta_{j_4}^{(i_4)}
\zeta_{j_5}^{(i_5)}-
$$
$$
-
{\bf 1}_{\{j_1=j_4\}}
{\bf 1}_{\{i_1=i_4\}}
\zeta_{j_2}^{(i_2)}
\zeta_{j_3}^{(i_3)}
\zeta_{j_5}^{(i_5)}-
{\bf 1}_{\{j_1=j_5\}}
{\bf 1}_{\{i_1=i_5\}}
\zeta_{j_2}^{(i_2)}
\zeta_{j_3}^{(i_3)}
\zeta_{j_4}^{(i_4)}-
$$
$$
-
{\bf 1}_{\{j_2=j_3\}}
{\bf 1}_{\{i_2=i_3\}}
\zeta_{j_1}^{(i_1)}
\zeta_{j_4}^{(i_4)}
\zeta_{j_5}^{(i_5)}-
{\bf 1}_{\{j_2=j_4\}}
{\bf 1}_{\{i_2=i_4\}}
\zeta_{j_1}^{(i_1)}
\zeta_{j_3}^{(i_3)}
\zeta_{j_5}^{(i_5)}-
$$
$$
-
{\bf 1}_{\{j_2=j_5\}}
{\bf 1}_{\{i_2=i_5\}}
\zeta_{j_1}^{(i_1)}
\zeta_{j_3}^{(i_3)}
\zeta_{j_4}^{(i_4)}-{\bf 1}_{\{j_3=j_4\}}
{\bf 1}_{\{i_3=i_4\}}
\zeta_{j_1}^{(i_1)}
\zeta_{j_2}^{(i_2)}
\zeta_{j_5}^{(i_5)}-
$$
$$
-
{\bf 1}_{\{j_3=j_5\}}
{\bf 1}_{\{i_3=i_5\}}
\zeta_{j_1}^{(i_1)}
\zeta_{j_2}^{(i_2)}
\zeta_{j_4}^{(i_4)}-{\bf 1}_{\{j_4=j_5\}}
{\bf 1}_{\{i_4=i_5\}}
\zeta_{j_1}^{(i_1)}
\zeta_{j_2}^{(i_2)}
\zeta_{j_3}^{(i_3)}+
$$
$$
+
{\bf 1}_{\{j_1=j_2\}}
{\bf 1}_{\{i_1=i_2\}}
{\bf 1}_{\{j_3=j_4\}}
{\bf 1}_{\{i_3=i_4\}}\zeta_{j_5}^{(i_5)}+
{\bf 1}_{\{j_1=j_2\}}
{\bf 1}_{\{i_1=i_2\}}
{\bf 1}_{\{j_3=j_5\}}
{\bf 1}_{\{i_3=i_5\}}\zeta_{j_4}^{(i_4)}+
$$
$$
+
{\bf 1}_{\{j_1=j_2\}}
{\bf 1}_{\{i_1=i_2\}}
{\bf 1}_{\{j_4=j_5\}}
{\bf 1}_{\{i_4=i_5\}}\zeta_{j_3}^{(i_3)}+
{\bf 1}_{\{j_1=j_3\}}
{\bf 1}_{\{i_1=i_3\}}
{\bf 1}_{\{j_2=j_4\}}
{\bf 1}_{\{i_2=i_4\}}\zeta_{j_5}^{(i_5)}+
$$
$$
+
{\bf 1}_{\{j_1=j_3\}}
{\bf 1}_{\{i_1=i_3\}}
{\bf 1}_{\{j_2=j_5\}}
{\bf 1}_{\{i_2=i_5\}}\zeta_{j_4}^{(i_4)}+
{\bf 1}_{\{j_1=j_3\}}
{\bf 1}_{\{i_1=i_3\}}
{\bf 1}_{\{j_4=j_5\}}
{\bf 1}_{\{i_4=i_5\}}\zeta_{j_2}^{(i_2)}+
$$
$$
+
{\bf 1}_{\{j_1=j_4\}}
{\bf 1}_{\{i_1=i_4\}}
{\bf 1}_{\{j_2=j_3\}}
{\bf 1}_{\{i_2=i_3\}}\zeta_{j_5}^{(i_5)}+
{\bf 1}_{\{j_1=j_4\}}
{\bf 1}_{\{i_1=i_4\}}
{\bf 1}_{\{j_2=j_5\}}
{\bf 1}_{\{i_2=i_5\}}\zeta_{j_3}^{(i_3)}+
$$
$$
+
{\bf 1}_{\{j_1=j_4\}}
{\bf 1}_{\{i_1=i_4\}}
{\bf 1}_{\{j_3=j_5\}}
{\bf 1}_{\{i_3=i_5\}}\zeta_{j_2}^{(i_2)}+
{\bf 1}_{\{j_1=j_5\}}
{\bf 1}_{\{i_1=i_5\}}
{\bf 1}_{\{j_2=j_3\}}
{\bf 1}_{\{i_2=i_3\}}\zeta_{j_4}^{(i_4)}+
$$
$$
+
{\bf 1}_{\{j_1=j_5\}}
{\bf 1}_{\{i_1=i_5\}}
{\bf 1}_{\{j_2=j_4\}}
{\bf 1}_{\{i_2=i_4\}}\zeta_{j_3}^{(i_3)}+
{\bf 1}_{\{j_1=j_5\}}
{\bf 1}_{\{i_1=i_5\}}
{\bf 1}_{\{j_3=j_4\}}
{\bf 1}_{\{i_3=i_4\}}\zeta_{j_2}^{(i_2)}+
$$
$$
+
{\bf 1}_{\{j_2=j_3\}}
{\bf 1}_{\{i_2=i_3\}}
{\bf 1}_{\{j_4=j_5\}}
{\bf 1}_{\{i_4=i_5\}}\zeta_{j_1}^{(i_1)}+
{\bf 1}_{\{j_2=j_4\}}
{\bf 1}_{\{i_2=i_4\}}
{\bf 1}_{\{j_3=j_5\}}
{\bf 1}_{\{i_3=i_5\}}\zeta_{j_1}^{(i_1)}+
$$
\begin{equation}
\label{yyy10}
+\Biggl.
{\bf 1}_{\{j_2=j_5\ne 0\}}
{\bf 1}_{\{i_2=i_5\}}
{\bf 1}_{\{j_3=j_4\ne 0\}}
{\bf 1}_{\{i_3=i_4\}}\zeta_{j_1}^{(i_1)}\Biggr),
\end{equation}

\vspace{7mm}
\noindent
where

\vspace{-1mm}
$$
C_{j_3j_2j_1}=\int\limits_{\tau_p}^{\tau_{p+1}}\phi_{j_3}(z)
\int\limits_{\tau_p}^{z}\phi_{j_2}(y)
\int\limits_{\tau_p}^{y}
\phi_{j_1}(x)dx dy dz=
$$

$$
=
\frac{\sqrt{(2j_1+1)(2j_2+1)(2j_3+1)}}{8}\Delta^{3/2}\bar
C_{j_3j_2j_1},
$$

\vspace{4mm}

$$
C_{j_4j_3j_2j_1}=\int\limits_{\tau_p}^{\tau_{p+1}}\phi_{j_4}(u)
\int\limits_{\tau_p}^{u}\phi_{j_3}(z)
\int\limits_{\tau_p}^{z}\phi_{j_2}(y)
\int\limits_{\tau_p}^{y}
\phi_{j_1}(x)dx dy dz du=
$$

$$
=\frac{\sqrt{(2j_1+1)(2j_2+1)(2j_3+1)(2j_4+1)}}{16}\Delta^{2}\bar
C_{j_4j_3j_2j_1},
$$

\vspace{4mm}

$$
C_{j_2j_1}^{01}
=\int\limits_{\tau_p}^{\tau_{p+1}}(\tau_p-y)\phi_{j_2}(y)
\int\limits_{\tau_p}^{y}
\phi_{j_1}(x)dx dy =
$$

$$
=\frac{\sqrt{(2j_1+1)(2j_2+1)}}{8}\Delta^{2}\bar
C_{j_2j_1}^{01},
$$

\vspace{4mm}
$$
C_{j_2j_1}^{10}
=\int\limits_{\tau_p}^{\tau_{p+1}}\phi_{j_2}(y)
\int\limits_{\tau_p}^{y}
(\tau_p-x)\phi_{j_1}(x)dx dy =
$$

$$
=
\frac{\sqrt{(2j_1+1)(2j_2+1)}}{8}\Delta^{2}\bar
C_{j_2j_1}^{10},
$$

\vspace{4mm}

$$
C_{j_3j_2j_1}^{001}=\int\limits_{\tau_p}^{\tau_{p+1}}(\tau_p-z)\phi_{j_3}(z)
\int\limits_{\tau_p}^{z}\phi_{j_2}(y)
\int\limits_{\tau_p}^{y}
\phi_{j_1}(x)dx dy dz=
$$

$$
=
\frac{\sqrt{(2j_1+1)(2j_2+1)(2j_3+1)}}{16}\Delta^{5/2}\bar
C_{j_3j_2j_1}^{001},
$$

\vspace{4mm}

$$
C_{j_3j_2j_1}^{010}=\int\limits_{\tau_p}^{\tau_{p+1}}\phi_{j_3}(z)
\int\limits_{\tau_p}^{z}(\tau_p-y)\phi_{j_2}(y)
\int\limits_{\tau_p}^{y}
\phi_{j_1}(x)dx dy dz=
$$

$$
=
\frac{\sqrt{(2j_1+1)(2j_2+1)(2j_3+1)}}{16}\Delta^{5/2}\bar
C_{j_3j_2j_1}^{010},
$$

\vspace{4mm}

$$
C_{j_3j_2j_1}^{100}=\int\limits_{\tau_p}^{\tau_{p+1}}\phi_{j_3}(z)
\int\limits_{\tau_p}^{z}\phi_{j_2}(y)
\int\limits_{\tau_p}^{y}
(\tau_p-x)\phi_{j_1}(x)dx dy dz=
$$

$$
=
\frac{\sqrt{(2j_1+1)(2j_2+1)(2j_3+1)}}{16}\Delta^{5/2}\bar
C_{j_3j_2j_1}^{100},
$$

\vspace{4mm}

$$
C_{j_5j_4 j_3 j_2 j_1}=
\int\limits_{\tau_p}^{\tau_{p+1}}\phi_{j_5}(v)
\int\limits_{\tau_p}^v\phi_{j_4}(u)
\int\limits_{\tau_p}^{u}
\phi_{j_3}(z)
\int\limits_{\tau_p}^{z}\phi_{j_2}(y)\int\limits_{\tau_p}^{y}\phi_{j_1}(x)
dxdydzdudv=
$$

$$
=\frac{\sqrt{(2j_1+1)(2j_2+1)(2j_3+1)(2j_4+1)(2j_5+1)}}{32}\Delta^{5/2}\bar
C_{j_5j_4 j_3 j_2 j_1},
$$

\vspace{5mm}
\noindent
where

\vspace{-2mm}
$$
\bar C_{j_3j_2j_1}=
\int\limits_{-1}^{1}P_{j_3}(z)
\int\limits_{-1}^{z}P_{j_2}(y)
\int\limits_{-1}^{y}
P_{j_1}(x)dx dy dz,
$$

\vspace{1mm}
$$
\bar C_{j_4j_3j_2j_1}=
\int\limits_{-1}^{1}P_{j_4}(u)
\int\limits_{-1}^{u}P_{j_3}(z)
\int\limits_{-1}^{z}P_{j_2}(y)
\int\limits_{-1}^{y}
P_{j_1}(x)dx dy dz,
$$

\vspace{1mm}

$$
\bar C_{j_2j_1}^{01}=-\int\limits_{-1}^{1}(1+y)P_{j_2}(y)
\int\limits_{-1}^{y}
P_{j_1}(x)dx dy,
$$

\vspace{1mm}
$$
\bar C_{j_2j_1}^{10}=-\int\limits_{-1}^{1}P_{j_2}(y)
\int\limits_{-1}^{y}
(1+x)P_{j_1}(x)dx dy,
$$

\vspace{1mm}

$$
\bar C_{j_3j_2j_1}^{100}=-
\int\limits_{-1}^{1}P_{j_3}(z)
\int\limits_{-1}^{z}P_{j_2}(y)
\int\limits_{-1}^{y}
P_{j_1}(x)(x+1)dx dy dz,
$$

\vspace{1mm}
$$
\bar C_{j_3j_2j_1}^{010}=-
\int\limits_{-1}^{1}P_{j_3}(z)
\int\limits_{-1}^{z}P_{j_2}(y)(y+1)
\int\limits_{-1}^{y}
P_{j_1}(x)dx dy dz,
$$

\vspace{1mm}

$$
\bar C_{j_3j_2j_1}^{001}=-
\int\limits_{-1}^{1}P_{j_3}(z)(z+1)
\int\limits_{-1}^{z}P_{j_2}(y)
\int\limits_{-1}^{y}
P_{j_1}(x)dx dy dz,
$$

\vspace{1mm}

$$
\bar C_{j_5j_4 j_3 j_2 j_1}=
\int\limits_{-1}^{1}P_{j_5}(v)
\int\limits_{-1}^{v}P_{j_4}(u)
\int\limits_{-1}^{u}P_{j_3}(z)
\int\limits_{-1}^{z}P_{j_2}(y)
\int\limits_{-1}^{y}
P_{j_1}(x)dx dy dz du dv,
$$

\vspace{6mm}
\noindent
where $P_i(x)$ $(i=0, 1, 2,\ldots)$ is the Legendre polynomial and

\vspace{2mm}

$$
\phi_i(x)=\sqrt{\frac{2i+1}{\Delta}}P_i\left(\left(x-\tau_p-\frac{\Delta}{2}\right)
\frac{2}{\Delta}\right),\ \ \ i=0, 1, 2,\ldots 
$$

\vspace{5mm}

Let us consider the exact relations and some estimates
for the mean-square approximation errors of iterated Ito
stochastic integrals.

Using Theorem 3, we get 
\cite{2017}-\cite{2013}, \cite{arxiv-3}, \cite{arxiv-24}

\vspace{1mm}
\begin{equation}
\label{xxx1}
{\sf M}\left\{\left(I_{(00)\tau_{p+1},\tau_p}^{(i_1 i_2)}-
I_{(00)\tau_{p+1},\tau_p}^{(i_1 i_2)q}
\right)^2\right\}
=\frac{\Delta^2}{2}\Biggl(\frac{1}{2}-\sum_{i=1}^q
\frac{1}{4i^2-1}\Biggr)\ \ \ (i_1\ne i_2),
\end{equation}

\vspace{7mm}

$$
{\sf M}\left\{\left(I_{(10)\tau_{p+1},\tau_p}^{(i_1 i_2)}-
I_{(10)\tau_{p+1},\tau_p}^{(i_1 i_2)q}
\right)^2\right\}=
{\sf M}\left\{\left(I_{(01)\tau_{p+1},\tau_p}^{(i_1 i_2)}-
I_{(01)\tau_{p+1},\tau_p}^{(i_1 i_2)q}\right)^2\right\}=
$$

\begin{equation}
\label{xxx2}
=\frac{\Delta^4}{16}\Biggl(\frac{5}{9}-
2\sum_{i=2}^q\frac{1}{4i^2-1}-
\sum_{i=1}^q
\frac{1}{(2i-1)^2(2i+3)^2}
-\sum_{i=0}^q\frac{(i+2)^2+(i+1)^2}{(2i+1)(2i+5)(2i+3)^2}
\Biggr)\ \ \ (i_1\ne i_2),
\end{equation}

\vspace{7mm}

$$
{\sf M}\left\{\left(I_{(10)\tau_{p+1},\tau_p}^{(i_1 i_1)}-
I_{(10)\tau_{p+1},\tau_p}^{(i_1 i_1)q}
\right)^2\right\}=
{\sf M}\left\{\left(I_{(01)\tau_{p+1},\tau_p}^{(i_1 i_1)}-
I_{(01)\tau_{p+1},\tau_p}^{(i_1 i_1)q}\right)^2\right\}=
$$

\begin{equation}
\label{xxx3}
=\frac{\Delta^4}{16}\Biggl(\frac{1}{9}-
\sum_{i=0}^{q}
\frac{1}{(2i+1)(2i+5)(2i+3)^2}
-2\sum_{i=1}^{q}
\frac{1}{(2i-1)^2(2i+3)^2}\Biggr).
\end{equation}

\vspace{8mm}

Applying (\ref{qq1})--(\ref{883}), we obtain

\vspace{2mm}

$$
{\sf M}\left\{\left(
I_{(000)\tau_{p+1},\tau_p}^{(i_1i_2 i_3)}-
I_{(000)\tau_{p+1},\tau_p}^{(i_1i_2 i_3)q}\right)^2\right\}=
\frac{\Delta^{3}}{6}-\sum_{j_3,j_2,j_1=0}^{q}
C_{j_3j_2j_1}^2\ \ \ (i_1\ne i_2,\  i_1\ne i_3,\ i_2\ne i_3),
$$

\vspace{3mm}

$$
{\sf M}\left\{\left(
I_{(000)\tau_{p+1},\tau_p}^{(i_1i_2 i_3)}-
I_{(000)\tau_{p+1},\tau_p}^{(i_1i_2 i_3)q}\right)^2\right\}=
\frac{\Delta^{3}}{6}-\sum_{j_3,j_2,j_1=0}^{q}
C_{j_3j_2j_1}^2
-\sum_{j_3,j_2,j_1=0}^{q}
C_{j_2j_3j_1}C_{j_3j_2j_1}\ \ \ (i_1\ne i_2=i_3),
$$

\vspace{3mm}

$$
{\sf M}\left\{\left(
I_{(000)\tau_{p+1},\tau_p}^{(i_1i_2 i_3)}-
I_{(000)\tau_{p+1},\tau_p}^{(i_1i_2 i_3)q}\right)^2\right\}=
\frac{\Delta^{3}}{6}-\sum_{j_3,j_2,j_1=0}^{q}
C_{j_3j_2j_1}^2
-\sum_{j_3,j_2,j_1=0}^{q}
C_{j_3j_2j_1}C_{j_1j_2j_3}\ \ \ (i_1=i_3\ne i_2),
$$

\vspace{3mm}

$$
{\sf M}\left\{\left(
I_{(000)\tau_{p+1},\tau_p}^{(i_1i_2 i_3)}-
I_{(000)\tau_{p+1},\tau_p}^{(i_1i_2 i_3)q}\right)^2\right\}=
\frac{\Delta^{3}}{6}-\sum_{j_3,j_2,j_1=0}^{q}
C_{j_3j_2j_1}^2
-\sum_{j_3,j_2,j_1=0}^{q}
C_{j_3j_1j_2}C_{j_3j_2j_1}\ \ \ (i_1=i_2\ne i_3).
$$

\vspace{8mm}

At the same time using the estimate (\ref{qq4}) 
for $i_1,\ldots,i_5=1,\ldots,m$, 
we have

\vspace{2mm}

$$
{\sf M}\left\{\left(
I_{(01)\tau_{p+1},\tau_p}^{(i_1i_2)}-
I_{(01)\tau_{p+1},\tau_p}^{(i_1i_2)q}\right)^2\right\}\le
2\Biggl(\frac{\Delta^{4}}{4}-\sum_{j_1,j_2=0}^{q}
\left(C_{j_2j_1}^{01}\right)^2\Biggr),
$$

\vspace{3mm}
$$
{\sf M}\left\{\left(
I_{(10)\tau_{p+1},\tau_p}^{(i_1i_2)}-
I_{(10)\tau_{p+1},\tau_p}^{(i_1i_2)q}\right)^2\right\}\le
2\Biggl(\frac{\Delta^{4}}{12}-\sum_{j_1,j_2=0}^{q}
\left(C_{j_2j_1}^{10}\right)^2\Biggr),
$$

\vspace{3mm}

\begin{equation}
\label{xxx4}
{\sf M}\left\{\left(
I_{(000)\tau_{p+1},\tau_p}^{(i_1i_2 i_3)}-
I_{(000)\tau_{p+1},\tau_p}^{(i_1i_2 i_3)q}\right)^2\right\}\le
6\Biggl(\frac{\Delta^{3}}{6}-\sum_{j_3,j_2,j_1=0}^{q}
C_{j_3j_2j_1}^2\Biggr),
\end{equation}

\vspace{3mm}

\begin{equation}
\label{xxx5}
{\sf M}\left\{\left(
I_{(0000)\tau_{p+1},\tau_p}^{(i_1i_2 i_3 i_4)}-
I_{(0000)\tau_{p+1},\tau_p}^{(i_1i_2 i_3 i_4)q}\right)^2\right\}\le
24\Biggl(\frac{\Delta^{4}}{24}-\sum_{j_1,j_2,j_3,j_4=0}^{q}
C_{j_4j_3j_2j_1}^2\Biggr),
\end{equation}

\vspace{3mm}

\begin{equation}
\label{xxx6}
{\sf M}\left\{\left(
I_{(100)\tau_{p+1},\tau_p}^{(i_1i_2 i_3)}-
I_{(100)\tau_{p+1},\tau_p}^{(i_1i_2 i_3)q}\right)^2\right\}\le
6\Biggl(\frac{\Delta^{5}}{60}-\sum_{j_1,j_2,j_3=0}^{q}
\left(C_{j_3j_2j_1}^{100}\right)^2\Biggr),
\end{equation}

\vspace{3mm}

\begin{equation}
\label{xxx7}
{\sf M}\left\{\left(
I_{(010)\tau_{p+1},\tau_p}^{(i_1i_2 i_3)}-
I_{(010)\tau_{p+1},\tau_p}^{(i_1i_2 i_3)q}\right)^2\right\}\le
6\Biggl(\frac{\Delta^{5}}{20}-\sum_{j_1,j_2,j_3=0}^{q}
\left(C_{j_3j_2j_1}^{010}\right)^2\Biggr),
\end{equation}
                          
\vspace{3mm}

\begin{equation}
\label{xxx8}
{\sf M}\left\{\left(
I_{(001)\tau_{p+1},\tau_p}^{(i_1i_2 i_3)}-
I_{(001)\tau_{p+1},\tau_p}^{(i_1i_2 i_3)q}\right)^2\right\}\le
6\Biggl(\frac{\Delta^5}{10}-\sum_{j_1,j_2,j_3=0}^{q}
\left(C_{j_3j_2j_1}^{001}\right)^2\Biggr),
\end{equation}

\vspace{3mm}

\begin{equation}
\label{xxx9}
{\sf M}\left\{\left(
I_{(00000)\tau_{p+1},\tau_p}^{(i_1 i_2 i_3 i_4 i_5)}-
I_{(00000)\tau_{p+1},\tau_p}^{(i_1 i_2 i_3 i_4 i_5)q}\right)^2\right\}\le
120\left(\frac{\Delta^{5}}{120}-\sum_{j_1,j_2,j_3,j_4,j_5=0}^{q}
C_{j_5 i_4 i_3 i_2 j_1}^2\right).
\end{equation}

\vspace{5mm}

As was emphasized in \cite{18a},
\cite{arxiv-14}, \cite{2018a}-\cite{2018aa1}, the Fourier--Legendre coefficients 
$\bar C_{j_3j_2j_1}$ and $\bar C_{j_4j_3j_2j_1}$ 
(as well as the Fourier--Legendre coefficients
$\bar C_{j_3j_2j_1}^{001},$ $\bar C_{j_3j_2j_1}^{010},$ 
$\bar C_{j_3j_2j_1}^{100},$
$\bar C_{j_5j_4j_3j_2j_1}$) can be exactly calculated using symbolic 
transfor\-ma\-ti\-on packages like
Derive. The exact values of these Fourier--Legendre coefficients calculated 
in Derive were presented in
tabular form in the monographs \cite{2006}-\cite{2013}. 
Note that the mendioned Fourier--Legendre coefficients 
do not depend
on the integration step $\tau_{p+1}-\tau_p$ of the numerical method, 
which can 
be variable.

Recently, 
the database with 270,000 exactly
calculated Fourier--Legendre coefficients was described \cite{Mikh-1}.
This database was used in the software package,
which is written in the Python programming language
for the implementation of explicit one-step strong numerical methods
with orders 0.5, 1.0, 1.5, 2.0, 2.5, and 3.0 
of convergence for Ito SDEs with non-commutative noise.
The optimization of the mean-square approximation 
procedures for iterated Ito stochastic integrals
from these numerical schemes can be found in \cite{Mikh-2}.

Generally speaking, the minimum values $q$ that guarantee the
fulfillment of the condition 
(\ref{4.3})
for each of approximations (see above)
are different and abruptly decreasing with the growth of 
orders of smallness with respect to $\Delta$ of
approximations of iterated Ito 
stochastic integrals.

For pairwise different $i_1,\ldots,i_5=1,\ldots,m$ 
Theorem 3 gives

\vspace{1mm}
\begin{equation}
\label{800}
{\sf M}\left\{\left(
I_{(000)\tau_{p+1},\tau_p}^{(i_1i_2 i_3)}-
I_{(000)\tau_{p+1},\tau_p}^{(i_1i_2 i_3)6}\right)^2\right\}=
\frac{\Delta^{3}}{6}-\sum_{j_1,j_2,j_3=0}^{6}
C_{j_3j_2j_1}^2
\approx
0.01956000\Delta^3,
\end{equation}

\vspace{2mm}
\begin{equation}
\label{801}
{\sf M}\left\{\left(
I_{(100)\tau_{p+1},\tau_p}^{(i_1i_2 i_3)}-
I_{(100)\tau_{p+1},\tau_p}^{(i_1i_2 i_3)2}\right)^2\right\}=
\frac{\Delta^{5}}{60}-\sum_{j_1,j_2,j_3=0}^{2}
\left(C_{j_3j_2j_1}^{100}\right)^2\approx
0.00815429\Delta^5,
\end{equation}

\vspace{2mm}
\begin{equation}
\label{802}
{\sf M}\left\{\left(
I_{(010)\tau_{p+1},\tau_p}^{(i_1i_2 i_3)}-
I_{(010)\tau_{p+1},\tau_p}^{(i_1i_2 i_3)2}\right)^2\right\}=
\frac{\Delta^{5}}{20}-\sum_{j_1,j_2,j_3=0}^{2}
\left(C_{j_3j_2j_1}^{010}\right)^2\approx
0.01739030\Delta^5,
\end{equation}

\vspace{2mm}
\begin{equation}
\label{803}
{\sf M}\left\{\left(
I_{(001)\tau_{p+1},\tau_p}^{(i_1i_2 i_3)}-
I_{(001)\tau_{p+1},\tau_p}^{(i_1i_2 i_3)2}\right)^2\right\}=
\frac{\Delta^5}{10}-\sum_{j_1,j_2,j_3=0}^{2}
\left(C_{j_3j_2j_1}^{001}\right)^2
\approx 0.02528010\Delta^5,
\end{equation}

\vspace{2mm}
\begin{equation}
\label{804}
{\sf M}\left\{\left(
I_{(0000)\tau_{p+1},\tau_p}^{(i_1i_2i_3 i_4)}-
I_{(0000)\tau_{p+1},\tau_p}^{(i_1i_2i_3 i_4)2}\right)^2\right\}=
\frac{\Delta^{4}}{24}-\sum_{j_1,j_2,j_3,j_4=0}^{2}
C_{j_4j_3j_2j_1}^2\approx
0.02360840\Delta^4,
\end{equation}

\vspace{2mm}
\begin{equation}
\label{805}
{\sf M}\left\{\left(
I_{(00000)\tau_{p+1},\tau_p}^{(i_1i_2i_3i_4 i_5)}-
I_{(00000)\tau_{p+1},\tau_p}^{(i_1i_2i_3i_4 i_5)1}\right)^2\right\}=
\frac{\Delta^5}{120}-\sum_{j_1,j_2,j_3,j_4,j_5=0}^{1}
C_{j_5j_4j_3j_2j_1}^2\approx
0.00759105\Delta^5.
\end{equation}

\vspace{5mm}

Recall that the value $\Delta$ acts as the integration step of the numerical 
method (\ref{4.45}) for the Ito
SDE (\ref{1.5.2}), thereby being rather small. Hence, 
even for $q=6, 2,$ and $1$ the
mean-square approximation errors (\ref{800})--(\ref{805}) of the iterated
Ito stochastic integrals of multiplicities
3 to 5 are sufficiently small.
Note that in \cite{KlPl2}, \cite{Mi2}, \cite{Mi3}
the iterated stochastic integrals were 
approximated using the trigonometric
Fourier expansion of the multidimensional
Brownian bridge process 
and the mean-square
approximation error of the 
iterated
stochastic integrals were estimated by the value

\vspace{-1mm}
$$
\frac{C_1\Delta^2}{q},
$$

\vspace{2mm}
\noindent
where $C_1$ is a constant and $\Delta,$ $q$ have the same 
meaning as in (\ref{qqqq1}). 
Clearly, such an approach
is rougher than the one involving Theorem 3.

Note that the number $q$ must be the same for all
approximations of iterated stochastic integrals from the
considered collection in the approach from 
\cite{KlPl2}, \cite{Mi2}, \cite{Mi3}
while the numbers $q$ can be 
chosen different for different stochastic integrals
from the
considered collection in the method based on Theorems 1--3.

On the basis of 
the presented 
expansions (see above) of 
iterated Ito stochastic integrals we 
can see that increasing of multiplicities of these integrals 
or degree indexes of their weight functions 
leads
to increasing 
of orders of smallness with respect to $\Delta$ in the mean-square sense 
for iterated stochastic integrals. This leads to a sharp decrease  
of member 
quantities (the numbers $q$)
in expansions of iterated Ito stochastic 
integrals, which are required for achieving the acceptable accuracy
of approximation.

\vspace{5mm}

\section{Numerical Algorithm With the Order 2.5 of Strong Convergence}

\vspace{5mm}

In this section, we will write the formulas and recommendations on 
the numerical method (\ref{4.45})
with the order 2.5 of strong convergence as an algorithm.

Let the Fourier--Legendre coefficients 

\vspace{-1mm}
$$
\bar C_{j_3j_2j_1},\ \ \ \bar C_{j_4j_3j_2j_1},\ \ \
\bar C_{j_3j_2j_1}^{001},\ \ \ \bar C_{j_3j_2j_1}^{010},\ \ \
\bar C_{j_3j_2j_1}^{100},\ \ \
\bar C_{j_5j_4j_3j_2j_1}
$$

\vspace{3mm}
\noindent
be precalculated \cite{Mikh-1}.

\vspace{3mm}

{\bf Algorithm 1.}

\vspace{3mm}

{\bf Step 1.} Specify the initial parameters of the problem: 
the integration interval $[0, T]$, the integration
step $\Delta$ (e.g., the constant one
$\Delta=T/N,$ where $N>1$; a variable step 
is also admissible), the
initial condition ${\bf y}_0$, 
and the constant $C$ appearing in the condition (\ref{4.3}).

\vspace{3mm}

{\bf Step 2.} Let $p=0.$

\vspace{3mm}

{\bf Step 3.} Choose the minimum values 
$q$ under which the right-hand
sides of (\ref{xxx1})--(\ref{xxx9}) are not exceeding the right-hand side of 
the inequality (\ref{4.3}).

\vspace{3mm}

{\bf Step 4.} Generate a sequence of independent standard Gaussian random 
variables $\zeta_l^{(i)}$ $(l=0, 1,\ldots,q+2;\ i=1,\ldots,m).$ Here
the number $q$ is a maximum from the numbers $q$ chosen at Step 3.

\vspace{3mm}

{\bf Step 5.} Model the iterated Ito stochastic integrals

$$
I_{(0)\tau_{p+1},\tau_p}^{(i_1)},\
I_{(1)\tau_{p+1},\tau_p}^{(i_1)},\ I_{(2)\tau_{p+1},\tau_p}^{(i_1)},\
I_{(00)\tau_{p+1},\tau_p}^{(i_1 i_2)},\ 
I_{(10)\tau_{p+1},\tau_p}^{(i_1 i_2)},\
I_{(01)\tau_{p+1},\tau_p}^{(i_1 i_2)},\
I_{(000)\tau_{p+1},\tau_p}^{(i_1 i_2 i_3)},
$$

\vspace{1mm}
$$
I_{(100)\tau_{p+1},\tau_p}^{(i_1 i_2 i_3)},\
I_{(010)\tau_{p+1},\tau_p}^{(i_1 i_2 i_3)},\
I_{(001)\tau_{p+1},\tau_p}^{(i_1 i_2 i_3)},\
I_{(0000)\tau_{p+1},\tau_p}^{(i_1 i_2 i_3 i_4)},\
I_{(00000)\tau_{p+1},\tau_p}^{(i_1 i_2 i_3 i_4 i_5)}
$$

\vspace{5mm}
\noindent
using the formulas (\ref{yyy1aaa})--(\ref{yyy10})
with the values $q$ 
chosen at Step 3.

\vspace{3mm}

{\bf Step 6.} Find ${\bf y}_{p+1}$ by the formula (\ref{4.45}).

\vspace{3mm}

{\bf Step 7.} If $p<N-1$, then assign $p=p+1$ and go back to Step 4. 
Otherwise proceed to Step 8.

\vspace{3mm}

{\bf Step 8.} End of the Algorithm 1.

\vspace{5mm}

\section{Conclusions}

In this paper, the efficient mean-square approximation 
procedures for the iterated Ito stochastic
integrals of multiplicities 1 to 5 that are based on the multiple 
Fourier--Legendre series have been developed.
These results can be used for the implementation of 
the numerical method (\ref{4.45}) with the
strong order 2.5 of convergence
for the Ito SDEs with multidimensional non-commutative noises.

\end{document}